\documentstyle[11pt]{article}

\newcommand{\cC}{{\cal C}}
\newcommand{\cD}{{\cal D}}
\newcommand{\cE}{{\cal E}}

\newcommand{\cH}{{\cal H}}
\newcommand{\cI}{{\cal I}}

\newcommand{\cK}{{\cal K}}

\newcommand{\cM}{{\cal M}}
\newcommand{\cN}{{\cal N}}
\newcommand{\cO}{{\cal O}}

\newcommand{\cQ}{{\cal Q}}
\newcommand{\cR}{{\cal R}}

\newcounter{example}

\newcommand{\B}[1]{\begin{#1}}
\newcommand{\E}[1]{\end{#1}}

\newcommand{\aX}{\mbox{$X_{an}$}}
\newcommand{\eX}{\mbox{$X_{\acute{e}t}$}}

\newcommand{\sa}{\mbox{\scriptsize{$\aX$}}} 

\renewcommand{\P}{\mbox{$\bf P$}}     
\renewcommand{\H}{\mbox{$\bf H$}}     
\newcommand{\A}{\mbox{$\bf A$}}     
\newcommand{\C}{\mbox{$\bf C$}}     
\newcommand{\Z}{\mbox{$\bf Z$}}     

\newcommand{\La}{\mbox{$\Lambda$}}

\newcommand{\Pic}{{\rm Pic}\,}

\newcommand{\Sp}{{\rm Spec}\,}
\newcommand{\df}{\mbox{\,$\stackrel{\pp{\rm def}}{=}$}\,}

\newcommand{\rank}{{\rm rank}\,}
\newcommand{\by}[1]{\stackrel{#1}{\rightarrow}}
\newcommand{\longby}[1]{\stackrel{#1}{\longrightarrow}}

\newcommand{\tensor}{\otimes}
\newcommand{\into}{\hookrightarrow}

\newcommand{\ie}{{\it i.e.\/}\ }

\newcommand{\p}[1]{\mbox{$\scriptstyle {#1}$}}
\newcommand{\pp}[1]{\mbox{$\scriptscriptstyle {#1}$}}
\newcommand{\cat}[1]{\mbox{\large $\cal {#1}$}}

\newcommand{\veq}{\mbox{\large $\parallel$}}

\newcommand{\vcong}{\cong\downarrow\;\;}
\newcommand{\limdir}[1]{{\displaystyle{\mathop{\rm
lim}_{\buildrel\longrightarrow\over{#1}}}}\,}

\newcommand{\bul}{\mbox{\Large $\cdot $}}

\title{$\cH$-cohomologies versus algebraic
cycles}

\author{{\sc L. Barbieri-Viale}\thanks{Partially
supported by the GNSAGA of CNR}} \date{}

\begin{document}

\maketitle

\begin{center}
\parbox{2.7in}{\small {\it Contents}\\ 0.\,Introduction
1.\,Preliminaries 2.\,Invariance  3.\,Homotopy
and proto-decomposition  4.\,Cap products  5.\,Cup
products  6.\,Intersection theory 7.\,Chern
classes and blow-ups\\ {\it Appendix}\, A.\,Examples and
comments\\ {\it References}}
\end{center}

\setcounter{section}{-1}

\section{Introduction}

After Quillen's proof of the Gersten conjecture (see
\cite{Q}), for algebraic regular schemes, a natural
approach to the theory of algebraic cycles appears to be
by dealing with the ``formalism'' associated to (local)
higher $K$-theory, as it is manifestly expressed by the
work of Bloch and Gillet (cf. \cite{BL}, \cite{GIN}). As
a matter of fact a more general and flexbile setting has
been exploited by Bloch and Ogus (see \cite{BO}) by
axiomatic methods.\\
The aim of this paper is to going further with this
axiomatic method in order to obtain a ``global
intersection theory'' (in the Grothendieck sense
\cite{GI}) directly from a given ``cohomology theory''.
To this aim we will assume given a (twisted) cohomology
theory $H^*$ and we let consider the $\cH$-cohomology
functor
$$X\leadsto H^{\#}_{Zar}(X,\cH^*\p{(\cdot)})$$
where $\cH^*\p{(\cdot)}$ is the Zariski sheaf associated
to $H^*$. By dealing with a cup-product structure on $H^*$
we are granted of a product in $\cH$-cohomology; by
arguing with the cap-product structure we are able to
obtain a cap-product between algebraic cycles and
$\cH$-cohomology classes, for $Y$ and $Z$ closed
subschemes of $X$  ($\La\df H^0({\rm point})$) $$\p{\cap}
:C_{n}(Y;\La)\tensor H^p_Z(X,\cH^p\p{(p)}) \to
C_{n-p}(Y\cap Z;\La)$$  where $C_{*}(-;\La)$ is the
``$\cH$-homology theory'' given by the hypercohomology of
the complexes of $E^1$-terms of the niveau spectral
sequence.\\ If $X$ is smooth of pure dimension $d$, by
capping with the ``fundamental cycle'' $[X]\in C_d(X;\La)$
we have a ``Poincar\'e duality'' isomorphism
$$[X]\p{\cap}- : H^p_Z(X,\cH^p\p{(p)}) \cong
C_{d-p}(Z;\La)$$
Thus the ``$\cH$-cycle class'' $\eta (Z)\in
H^p_Z(X,\cH^p\p{(p)})$ is defined by $[X]\p{\cap}\eta (Z) =
[Z]$, for $i: Z\into X$ a closed subscheme of pure
codimension $p$ in $X$. By capping with the $\cH$-cycle
class we do obtain Gysin maps for algebraic cycles i.e.
maps $i^!:C_{n}(X;\La)\to C_{n-p}(Z;\La)$. Furthermore,
the $\cH$-cycle classes are compatible with the
intersection of cycles (when existing!) so that the
$\cH$-cohomology rings generalize the classical
``intersection rings'' obtained via rational or algebraic
equivalences (cf. \cite{GIN} for the $\cK$-cohomology).\\
The covariant property of the niveau spectral sequence
grant us of ``$\cH$-Gysin maps'' $$f_*:
H^{\#}_{f^{-1}(Z)}(Y,\cH^*\p{(\cdot)}) \to
H^{\#+\rho}_{Z}(X,\cH^{*+\rho}\p{(\cdot+\rho)})$$
associated with a proper morphism $f:Y\to X$ of relative
dimension $\rho$ between smooth schemes. The corresponding
projection formula holds. By the homotopy property of
$H^*$ we are obtaining homotopy and Dold-Thom
decomposition for $\cH$-cohomologies. By observing that
the canonical cycle map for line bundles $c\ell : \Pic (X)
\to H^2(X,\p{1})$ has always its image contained in the
subgroup of the locally trivial cohomology classes i.e.
$H^1(X,\cH^1\p{(1)})$ by the coniveau spectral sequence,
we are able to construct Chern classes in
$\cH$-cohomologies according with Gillet and Grothendieck
(see \cite{GIL} and \cite{GC}) $$c_{p,i} : K_i^Z(X) \to
H^{p-i}_Z(X,\cH^p\p{(p)})$$ where $Z$ is any closed subset
of $X$ smooth. These yield Riemann-Roch theorems and,
notably, Chern classes in $H^{2*}(-,\p{*})$ by composition
with the cycle map $H^*(-,\cH^*\p{(*)})\to
H^{2*}(-,\p{*})$ canonically induced by the coniveau
spectral sequence.\\ At last, an immediate application of
this setting is the ``blow-up formula''
$$H^p(X',\cH^q\p{(r)}) \cong H^p(X,\cH^q\p{(r)})\oplus
\bigoplus_{i=0}^{c-2}
H^{p-1-i}(Z,\cH^{q-1-i}\p{(r-1-i)})$$
where $X'$ is the blow-up of $X$ smooth, along a closed
smooth subset $Z$ of pure codimension $c$. Remarkably the
formula is obtained by no use of ``self-intersection''
nor ``formule-clef'' (used by the redundant arguments
made in \cite[Expos\'e VII]{SGA5} for \'etale cohomology or
Chow groups).\\
The paper is organized by adding structure to the assumed
Bloch-Ogus cohomology to proving the claimed results. The
common cohomologies (e.g. \'etale, de Rham or
Deligne-Beilinson cohomology) are examples for this setting
as explained in the Appendix. Some of these results are
already been used by the author for applications (see
\cite{BV1},\cite{BV2}); as our second main goal is
that this ``formalism'' can be used for applications to
birational geometry and algebraic cycles.\\ I would like
to thank H.Gillet and R.W.Thomason for some corroborating
conversations on these topics.

\section{Preliminaries}

Let $\cat{V}_k$ be the category of schemes of finite type
over a fixed ground field $k$; usually, an object of
$\cat{V}_k$ is called `algebraic scheme'. Let
$\cat{V}^2_k$ be the category whose objects are pairs
$(X,Z)$ where $X$ is an algebraic scheme and $Z$ is a
closed subscheme of $X$; morphisms in $\cat{V}^2_k$ are
fibre products in $\cat{V}_k\ .$ In the following we will
consider a contravariant functor $$(X,Z) \leadsto
H_Z^*(X,\cdot)$$ from  $\cat{V}^2_k$ to \Z -bigraded
abelian groups.\\ We need to assume, at least, that the
above functor gives rise to a `Poincar\'e duality theory
with supports' as setted out by Bloch and Ogus
\cite[1.1-1.3 and 7.1.2]{BO}; one can also consider
$K_i^Z(X)$ the relative Quillen $K$-theory \cite{Q} (see
\cite[Def. 2.13]{GIL}). For $X\in \cat{V}_k$ we denote
$H^*(X,\cdot)$ for $H_X^*(X,\cdot)\ .$

\subsection{Bloch-Ogus theory}\label{inter1}

 For the sake of notation we recall some facts by
\cite{BO}. Togheter with the cohomology theory $H^*(\
,\cdot)$ is given an homological functor $H_{\star}(\
,\dagger)$ covariant for proper morphisms in $\cat{V}_k$
and a pairing:  $$\p{\cap_{X,Z}}: H_l(X,\p{m})\tensor
H_Z^r(X,\p{s})\to H_{l-r}(Z,\p{m-s})$$
having a `projection formula'. For $f$ a proper map let
$f_!$ denote the induced map on homology.
It is also assumed the existence of a `fundamental class'
$\eta_X\in H_{2d}(X,\p{d})$ , $d = {\rm dim} X$ , such
that $f_!(\eta_X)=[K(X):K(Y)]\cdot \eta_Y$
if $f:X\to Y$ is proper and  ${\rm dim} X ={\rm dim} Y$
(cf.~\cite[7.1.2]{BO}). For $X$ smooth of dimension $d$
 $$\eta_X\p{\cap_{X,Z}}- :
H_Z^{2d-i}(X,\p{d-j})\by{\simeq} H_{i}(Z,\p{j})$$
is an isomorphism (`Poincar\'e duality') suitably
compatible with restrictions (cf.~\cite[1.4]{BO}). For
$Z\subseteq T\subseteq X$ such that $Z$ and $T$ are closed
in $X$ there is a long exact sequence
(see \cite[1.1.1]{BO})
 \B{equation}\label{loc}
\cdots \to H_Z^i(X,\cdot)\to H_T^i(X,\cdot)\to
H_{T-Z}^i(X-Z,\cdot)\to H_Z^{i+1}(X,\cdot)\to \cdots
\E{equation} suitably contravariant; moreover, for $X$
smooth and irreducible of dimension $d$ , the following
exact sequence ( $h+i=2d\ ,\ \dagger + \cdot =d$ ):
 \B{equation}\label{hloc} \cdots \to
H_h(Z,\dagger)\to H_h(T,\dagger)\to H_h(T-Z,\dagger)\to
H_{h-1}(Z,\dagger)\to \cdots \E{equation} is the
corresponding Poincar\'e dual of the above (cf. \cite[6.1
k)]{JA}).

\subsection{Gersten or arithmetic resolution}

Let $Z^p(X)= \{Z\subset X \mbox{: closed of}\ {\rm
codim}_XZ\geq p\}$, ordered by inclusion, and let define
$$H^i_{Z^p(X)}(X,\cdot) \df \limdir{\pp{Z\in Z^p(X)}}
H_Z^i(X,\cdot) $$ and for $x\in X$  $$H^i(x,\cdot) \df
\limdir{\pp{ U open \subset \overline{\{ x\} }}}
H^i(U,\cdot)$$ Taking the direct limit of the exact
sequences (\ref{loc}) (over the pairs $Z\subseteq T$ with
$Z\in Z^{p+1}(X)$ and $T\in Z^p(X)$) and using
`local purity' (cf. \cite[Prop.3.9]{BO}) on $X$ smooth
over $k$ perfect, one obtains long exact sequences
\B{equation}\label{limloc}
 H^i_{Z^{p+1}(X)}(X,\p{j})\to
H^i_{Z^p(X)}(X,\p{j})\to \coprod_{x\in X^p}^{}
H^{i-2p}(x,\p{j-p}) \to H^{i+1}_{Z^{p+1}(X)}(X,\p{j})
\E{equation}
where $X^p$ is the set of points whose closure
has codimension $p$ in $X$. Furthermore, if $f:X \to Y$ is
a flat morphism and $Z\in Z^p(Y)$ then we have $
f^{-1}(Z)\in Z^p(X)$; thus the sequence in (\ref{limloc})
yields a sequence of presheaves for the Zariski topology.
Let $H^*_{Z^p,X}\p{(\cdot)}$ denote the Zariski presheaf
$U \leadsto  H^*_{Z^p(U)}(U,\cdot)$ on $X$ and let $a:
\cat{P}(X_{Zar}) \to \cat{S}(X_{Zar})$ be the associated
sheaf exact functor. Denote $$aH^*_{Z^p,X}\p{(\cdot)}\df
\cH^*_{Z^p,X}\p{(\cdot)}$$ The presheaf
$H^*_{Z^0,X}\p{(\cdot)}$ is just the functor
$H^*\p{(\cdot)}$ on $X_{Zar}$ and so one has
$\cH^*_{Z^0,X}\p{(\cdot)} = \cH^*_X\p{(\cdot)}\ .$ One of
the main results of \cite{BO} is in proving the
vanishing of the map $$\cH^*_{Z^{p+1},X}\p{(\cdot)} \to
\cH^*_{Z^p,X}\p{(\cdot)}$$ for all $p \geq 0$ . From this
vanishing, sheafifying the sequence (\ref{limloc}),
one has the following exact sequences of sheaves on $X$
smooth over $k$ perfect:
\B{equation}\label{shortloc} 0\to
\cH^i_{Z^{p},X}\p{(j)}\to \coprod_{x\in X^p}^{} i_x
H^{i-2p}(x,\p{j-p}) \to  \cH^{i+1}_{Z^{p+1},X}\p{(j)}\to 0
\E{equation}
where: for $A$ an abelian group and $x\in X$ we let
$i_xA$ denote the constant sheaf $A$ on $\overline{\{
x\}}$ extended by zero to all $X$. Patching toghether
the above short exact sequences we do get a resolution of
the sheaf $\cH^i_X\p{(j)}$ (`arithmetic resolution'
in \cite[Theor.4.2]{BO}):
$$
0\to \cH^i_X\p{(j)}\to \coprod_{x\in X^0}^{}
i_x H^{i}(x,\p{j}) \to \coprod_{x\in X^1}^{}
i_x H^{i-1}(x,\p{j-1}) \to \cdots
$$

\B{rmk}
The assumption of $k$ perfect is unnecessary (cf.
\ref{app}.1) if $H^*(\ ,\cdot)$ is the \'etale theory
(namely $H^i(X,\p{j}) \df H^i(\eX,\mu_{\nu}^{\tensor j})$
where $\mu_{\nu}$ is the \'etale sheaf of $\nu^{\rm th}$
root of unity and $\nu$ is any positive integer prime to
char($k$) ).
 \E{rmk}

\subsection{Quillen $K$-theory}\label{inter2}

Let $X\leadsto K_p(X)$ be the Quillen $K$-functor
associated with the exact category of vector bundles on
any scheme $X$ (see \cite{Q}). For a fixed $X$ we let
$\cK_{p}(\cO_X)$ be the associated Zariski sheaf on $X$.
For any noetherian separated scheme $X$ we have a complex
of flasque sheaves  (`Gersten's complex')
 $$\cI^{\bul}_{q,X}: \coprod_{x\in X^0}^{} i_x K_{q}(k(x))
\to \coprod_{x\in X^1}^{} i_x K_{q-1}(k(x)) \to \cdots $$
conjecturally exact if $X$ is regular, being proved exact
by Quillen \cite{Q} if $X$ is regular and essentially of
finite type over a field. By tensoring coherent modules
with locally free sheaves and sheafifying one has a
pairing of complexes of sheaves for $p,q\geq 0$ (see
\cite[p.276-277]{GIL}):
$$\cap : \cK_{p}(\cO_X)\tensor
\cI^{\bul}_{q,X}\to \cI^{\bul}_{p+q,X}$$
By capping with $1_X$ = the `fundamental class'
i.e. the identity section of the constant sheaf $\Z
\cong \cI^{\bul}_{0,X}$ we get an augmentation
$$\cap 1_X : \cK_{p}(\cO_X) \to \cI^{\bul}_{p,X}$$
which is a quasi-isomorphism if $X$ is regular essentially
of finite type over a field (conjecturally for all
regular schemes.) If $f:X\to Y$ is a proper morphism
between biequidimensional schemes and  $r={\rm dim}\,Y -
{\rm dim}\,X$ then there is an induced map of complexes
 $f_{!}: f_*\cI^{\bul}_{q,X} \to \cI^{\bul}_{q+r,Y}[r]$
(which takes the elements of $K_*(k(x))$ to
$K_*(k(f(x)))$ if dim $\bar{x}$ = dim $\overline{f(x)}$ and
takes them to zero otherwise) and a commutative diagram
(see the `projection formula'  \cite[p.411]{GIN}):

\B{equation}\B{array}{c}\label{kpr}
\hspace{20pt}f_*\cK_{p}(\cO_X)\tensor
f_*\cI^{\bul}_{q,X}\\
\p{f^{\natural}\tensor id}\nearrow
\ \hspace{40pt}\ \searrow \p{\cap_{f_*}} \\
\cK_{p}(\cO_Y)\tensor f_*\cI^{\bul}_{q,X}\hspace{60pt}
f_*\cI^{\bul}_{p+q,X}\\ \p{id\tensor f_{!}}\downarrow \
\hspace{105pt}  \ \downarrow \p{f_{!}} \\
\cK_{p}(\cO_Y)\tensor \cI^{\bul}_{q+r,Y}[r]  \hspace{10pt}
\longby{\cap}  \hspace{10pt} \cI^{\bul}_{p+q+r,Y}[r]
\E{array}
\E{equation}

Thus the formula: $f_!(f^{\natural}(\tau)
\cap_{f_*}\sigma) = \tau\cap f_!(\sigma )  $
for all sections $\tau \p{\tensor} \sigma$ of the
complex of sheaves $\cK_{p}(\cO_Y)\tensor
f_*\cI^{\bul}_{q,X}$\ .

\section{Invariance}

Let $X\in \cat{V}_k$ be an algebraic scheme; in the
following we will assume $X$ smooth and the ground field
$k$ perfect. We moreover assume given a cohomological
functor $H_Z^*(X,\cdot)$ satisfying the list of axioms
\cite[1.1--1.3]{BO} and the assumption \cite[7.1.2]{BO}.
We are going to consider a morphism $f:X\to Y$
from $X$ as above to $Y\in \cat{V}_k$ tacitly assuming $Y$
to be smooth.\\ Let X be in $\cat{V}_k$ and let
$\cH^i_X\p{(j)}$ (resp. $\cK_{i}(\cO_X)$) denote the sheaf
on $X$, for the Zariski topology, associated to the
presheaf $U \leadsto H^i(U, \p{j})$ (resp. $U \leadsto
K_i(\mbox{vector bundles on}\, U)$). If $f:X\to Y$ is any
morphism in $\cat{V}_k$ then there are maps
$f^{\sharp}:  \cH^i_Y\p{(j)} \to f_*\cH^i_X\p{(j)}$ (resp.
$f^{\natural}:  \cK_{i}(\cO_Y) \to f_*\cK_{i}(\cO_X)$).

\B{teor} {\rm (Invariance)}
Let $f:X\to Y$ be a proper birational morphism in
$\cat{V}_k$. For $X$ and $Y$ smooths over $k$
perfect then $f^{\sharp}$  yields
the isomorphism $$\cH^i_Y\p{(j)}\cong f_*\cH^i_X\p{(j)}$$
for all integers $i$ and $j$.
For $X$ and $Y$ regular algebraic schemes $f^{\natural}$
induces the isomorphism
$$ \cK_{i}(\cO_Y)\cong f_*\cK_{i}(\cO_X)$$
for all $i\geq 0$.\\
Hence there are isomorphisms
$$H^0(X,\cH^i_X\p{(j)})\cong H^0(Y,\cH^i_Y\p{(j)})$$
and
$$ H^0(X,\cK_{i}(\cO_X))\cong H^0(Y,\cK_{i}(\cO_Y))$$
\E{teor}

\B{rmk} By the recent results of M. Spivakovsky
the problem of `the elimination of
points of indeterminacy' appears to be solved also in
positive characteristic. Thus, by the Theorem~1, the groups
$H^0(X,\cH^{*}_X\p{(\cdot)})$ and $H^0(X,\cK_{*}(\cO_X))$
are birational invariants of $X$ smooth and proper (cf.
\cite{CO}).
\E{rmk}

The proof of the Theorem~1 is quite natural and easy after
a sheaf form of the projection formula (cf.
Lemma~\ref{chiave} and (\ref{hpr})). We will
first give an explicit description of the map $f^{\sharp}$.

\subsection{Functoriality}\label{rip1}

If we are given a morphism $f:X\to Y$ then,
for any  open $V$ subset of $Y$, there is a homomorphism
$H^*(V,\cdot) \to H^*(f^{-1}(V),\cdot)$ induced by $f$
simply because $H^*(\ ,\cdot)$ is a contravariant
functor; thus, with the notation previously introduced, we
get indeed a map $$f^H:H^*_Y\p{(\cdot)} \to
f_*H^*_X\p{(\cdot)}$$ of presheaves on $Y$. Moreover there
is a canonical map $f_*H^*_X\p{(\cdot)} \to
f_*\cH^*_X\p{(\cdot)}$ induced by sheafification and
direct image. Hence, taking the associated sheaves, one get
\B{eqnarray}\label{map}
\cH^*_Y\p{(\cdot)} & \by{af^H} & af_*H^*_X\p{(\cdot)}
\nonumber \\
     &  & \  \downarrow  \p{f_a} \\
     & & f_*\cH^*_X\p{(\cdot)} \nonumber
\E{eqnarray}
Thus the map $f^{\sharp}: \cH^*_Y\p{(\cdot)} \to
f_*\cH^*_X\p{(\cdot)}$ is defined to be the composite of
$af^H$ and $f_a$ as above.

\subsection{Key Lemma}\label{svol1}

In the following, till the end of this subsection, we will
let $f:X\to Y$ be a proper birational morphism between
smooth algebraic schemes over a perfect field. Our goal is
to prove that $f^{\sharp }$ is an isomorphism of
sheaves.\\[1pt]

{\it Step 1}.\, We can reduce to proving Theorem~1 for
irreducible schemes because if not then, from smoothness,
the irreducible components coincide with the connected
components and $f$ maps components to components;
hence, if $X_0$ and $Y_0$ are components such that
$f:X_0\to Y_0$ and $y\in Y_0$ then there are
isomorphisms on the stalk $(\cH^*_Y\p{(\cdot)})_y\cong
(\cH^*_{Y_0}\p{(\cdot)})_y$ and $(f_*\cH^*_X\p{(\cdot)})_y\cong
(f_*\cH^*_{X_0}\p{(\cdot)})_y$
If $X$ is irreducible and $K(X)$ is the function field of
$X$ then we denote  $$H^*(K(X),\cdot) \df \limdir{U
open\subset X} H^*(U,\cdot).$$ $H^*(K(X),\cdot)$ is
canonically contravariant and birationally invariant.
Hence $f$ induces an isomorphism $H^*(K(Y),\cdot)\cong
H^*(K(X),\cdot).$ So, we moreover assume $X$ and $Y$
irreducibles.\\[1pt]

{\it Step 2}.\, Case $\cH^0\p{(\cdot)}.$ Assume that the
cohomology theory is concentrated in positive degrees \ie
$H^i\p{(\cdot)}=0$ if $i<0$; hence the arithmetic
resolution yields an isomorphism $\cH^0_X\p{(\cdot)}\cong
i_XH^0(K(X),\cdot)$ = the constant sheaf $H^0(K(X),\cdot)$
on $X$. The same holds on $Y$. Since $f$ has connected
fibres (`Zariski main theorem') then
$\cH^0_Y\p{(\cdot)}\cong f_*\cH^0_X\p{(\cdot)}\ .$ The non
bounded case is considered below.\\[1pt]

{\it Step 3}.\, So, associated to $f:X\to Y$, by
(\ref{shortloc}) and (\ref{map}), we can construct a
diagram  \B{equation} \label{diagramma chiave}
\begin{array}{ccccccccc}
{0}&{\to}&{\cH^i_{Y}\p{(j)}}&{\to}&{i_YH^i(K(Y),\p{j})} &
{\to}& {\cH^{i+1}_{Z^{1},Y}\p{(j)}} &{\to}&{0} \\
   & &{\p{f^{\sharp}}\downarrow \ \ } & &{\vcong} & &
{\downarrow \ \ \p{f^{\sharp}_{Z^1}}} & &\\
{0}&{\to}&{f_*\cH^i_{X}\p{(j)}}&{
\to}&{f_*(i_XH^i(K(X),\p{j}))}&{\to}&
{f_*(\cH^{i+1}_{Z^{1},X}\p{(j)})}& &
\end{array}
\E{equation}
where the right most vertical arrow (it will
be seen explicitly below) is defined by commutativity of
the left hand square. (Note: because $f$ has connected
fibres then the middle vertical map is an isomorphism. The
commutativity is straightforward.)\\
{}From the above diagram one can see that $f^{\sharp}:
\cH^*_Y\p{(\cdot)} \to f_*\cH^*_X\p{(\cdot)}$ is
injective. Because of $f$ proper, and the arithmetic
resolution is covariant for proper maps, we do aim to get
the following commuative diagram
 \B{displaymath}\B{array}{ccccccc}  0
&\to & \cH^{i}_{Y}\p{(j)} & \to &{i_YH^i(K(Y),\p{j})
}&\to &{\displaystyle \coprod_{y\in Y^1}^{} i_y
H^{i-1}(y,\p{j-1})}\\
 & &\p{f_{\sharp}} \uparrow \ \ & & \cong \uparrow \ \ & &
\uparrow \\ 0 & \to & f_*(\cH^{i}_{X}\p{(j)}) &\to
&  f_*(i_XH^i(K(X),\p{j})) &\to & f_*({\displaystyle
\coprod_{x\in X^1}^{} i_x H^{i-1}(x,\p{j-1})})
\E{array}\E{displaymath}
where $f_{\sharp}$ is an injection.
The Theorem~1 is obtained by proving:
$f_{\sharp}\p{\circ }f^{\sharp}=id$ as a
consequence of the projection formula. Indeed we have:

\B{lemma}\label{chiave}
Let $f:X\to Y$ be a proper birational morphism
between irreducible algebraic smooth schemes. Then there
are maps of sheaves ($k=0,1$):
$$f^{\sharp}_{Z^k}:\cH^{i}_{Z^{k},Y}\p{(j)} \to
f_*(\cH^{i}_{Z^{k},X}\p{(j)})$$  and
$$f_{\sharp}^{Z^k}:f_*(\cH^{i}_{Z^{k},X}\p{(j)}) \to
\cH^{i}_{Z^{k},Y}\p{(j)}$$ such that
$$f_{\sharp}^{Z^k}\p{\circ }f^{\sharp}_{Z^k}=id$$
\E{lemma}
(Remind: $f_{\sharp}^{Z^0} =f_{\sharp}$ and
$f^{\sharp}_{Z^0}=f^{\sharp}$.)\\[2pt]
\B{proof}
We will follow the framework given by
Grothendieck in \cite[III.9.2]{GR}.\\ Note that $f(X)$ is
closed and dense in $Y$ irreducible: $f(X)=Y$. For all
$Z\in Z^k(Y)$ so that $f^{-1}(Z)\in Z^k(X)$ ($k=0,1$) we
have a map $H^{i}_{Z}(Y,\cdot)\to
H^{i}_{f^{-1}(Z)}(X,\cdot)$ and since $f$ is a proper
morphism between smooth schemes we have also maps
$H^{i}_{f^{-1}(Z)}(X,\cdot)\to H^{i}_{Z}(Y,\cdot)$ for
all $Z\in Z^k(Y)$. We then have: \B{slemma}\label{prfor}
The composition: $$H^{i}_{Z}(Y,\cdot) \by{f^{\star}}
H^{i}_{f^{-1}(Z)}(X,\cdot) \by{f_{\star}}
H^{i}_{Z}(Y,\cdot)$$ is the identity.
\E{slemma}
\B{proof}
Let $H_*(\ ,\dagger)$ denote the `twin' homology theory and
consider the pairing: $$\p{\cap_{Y,Z}}:
H_l(Y,\p{m})\tensor H_Z^r(Y,\p{s})\to H_{l-r}(Z,\p{m-s})$$
Denote $f_{!}: H_*(f^{-1}(T),\dagger) \to  H_*(T,\dagger)$
the homomorphisms induced, by covariance, from the proper
maps  $f^{-1}(T) \to T $ for every closed subset $T$ of
$Y$. Let $f^{\star}: H^*_{Z}(Y,\cdot) \to
H^{*}_{f^{-1}(Z)}(X,\cdot)$ be the map given by
contravariancy. Because of \cite[Axiom 1.3.3]{BO}  we have
the projection formula:
$$f_!(x\p{\cap_{X,f^{-1}(Z)}}f^{\star}(y))=f_!(x)\p{
\cap_{Y,Z}}y$$ for every $x\in H_l(X,m)$ and $y\in
H_Z^r(Y,s)$. Let $\eta_X$ denote the fundamental class
in $H_{2d}(X,d)$ where $d= {\rm dim} X$. Because of
\cite[7.1.2]{BO} and $f$ proper birational we get:
$f_!(\eta_X)=\eta_Y$. Thus the projection formula yields
the equation:
$$f_!(\eta_X\p{\cap_{X,f^{-1}(Z)}}f^{\star}(y))=
\eta_Y\p{\cap_{Y,Z}}y$$  By Poincar\'e duality
\cite[1.3.5]{BO} the cap product  with the fundamental
class is an isomorphism; we define $$f_{\star}(z) \df
(\eta_Y\p{\cap_{Y,Z}}-)^{-1}\p{\circ}
f_!(\eta_X\p{\cap_{X,f^{-1}(Z)}}z)$$ for all $z\in
H^{*}_{f^{-1}(Z)}(X,\cdot)$. Thus: $f_{\star}\p{\circ}
f^{\star}=1.$
\E{proof}

Taking the direct limit of the concerned maps over
$Z\in Z^k(Y)$ (note: because $f$ is closed the
direct system $\{ f^{-1}(Z)\, :\, Z\in Z^1(Y)\}$ is
cofinal in $Z^1(X)$) we have that the composition
$$H^{i}_{Z^k(Y)}(Y,\cdot)\by{f^{\star}_{Z^k}}
H^{i}_{Z^k(X)}(X,\cdot)
\by{f_{\star}^{Z^k}}H^{i}_{Z^k(Y)}(Y,\cdot)$$
is the identity as a consequence of the
Sublemma~\ref{prfor} and limit arguments (the
compatibilities are given by
\cite[1.1.2~and~1.2.4]{BO}).\\  Because
of~\cite[1.2.2~and~1.4]{BO} the maps $f_{\star}^{Z^k}$ are
natural trasformations of Zariski presheaves
$H^*_{Z^k,Y}(\cdot) \to f_*H^*_{Z^k,X}(\cdot)$ on $Y$.
Thus, taking the associated sheaves, we have that:
$$\cH^{i}_{Z^{k},Y}\p{(j)} \by{af^{\star}_{Z^k}}
af_*(H^{i}_{Z^{k},X}\p{(j)}) \by{af_{\star}^{Z^k}}
\cH^{i}_{Z^{k},Y}\p{(j)} $$ is the identity. Now it
sufficies to make up a commutative diagram as follows
\B{equation}\B{array}{ccc}\label{bravo}
\cH^{i}_{Z^{k},Y}\p{(j)}\ \by{af^{\star}_{Z^k}} &
af_*(H^{i}_{Z^{k},X}\p{(j)}) & \by{af_{\star}^{Z^k}} \
\cH^{i}_{Z^{k},Y}\p{(j)} \\ \p{f^{\sharp}_{Z^k}}\searrow
& \  \downarrow  \p{f_a^{Z^k}}&  \nearrow
\p{f_{\sharp}^{Z^k}} \\
      & f_*(\cH^{i}_{Z^{k},X}\p{(j)}) &
\E{array}\E{equation}
{}From (\ref{bravo}) we then have:
$$f_{\sharp}^{Z^k}\p{\circ}f^{\sharp}_{Z^k}=
f_{\sharp}^{Z^k}\p{\circ}\,f_a^{Z^k}\p{\circ}\,
af^{\star}_{Z^k}=
af_{\star}^{Z^k}\p{\circ}\,af^{\star}_{Z^k}=id$$  as
claimed. Indeed $f^{\sharp}_{Z^k}$ is simply defined by
composition; since $f$ is proper, dim$X$ = dim$Y$ and the
arithmetic resolution is covariant for proper maps:
$f_{\sharp}^{Z^k}$ is obtained, e.g.  $f_{\sharp}^{Z^1}$
from the commutativity and the exactness of the following:
\B{displaymath}\B{array}{ccccccc}  0 &\to &
\cH^{i+1}_{Z^{1},Y}\p{(j)} & \to &{\displaystyle
\coprod_{y\in Y^1}^{} i_y H^{i-1}(y,\p{j-1}) }&\to
&{\displaystyle \coprod_{y\in Y^2}^{} i_y
H^{i-2}(y,\p{j-2})}\\  & &  & & \uparrow & & \uparrow \\ 0
& \to & f_*(\cH^{i+1}_{Z^{1},X}\p{(j)}) &\to &
f_*({\displaystyle \coprod_{x\in X^1}^{} i_x
H^{i-1}(x,\p{j-1})}) &\to & f_*({\displaystyle
\coprod_{x\in X^2}^{} i_x H^{i-2}(x,\p{j-2})})
\E{array}\E{displaymath}  The proof of the
Lemma~\ref{chiave} is complete. \E{proof}

\subsection{Proof of the Invariance Theorem}\label{H-coda}

To summarize the proof: if $k=0,1$ and
$Z\in Z^k(Y)$ then $f^{-1}(Z)\in Z^k(X)$; we
have a splitting between long exact sequences
(cf.~(\ref{loc}))    $$
\begin{array}{ccccccccc}
{\cdots}&{\to}&{H^i(Y,\cdot)}&{\to}&{H^i(Y- Z,\cdot)} &
{\to}& {H^{i+1}_{Z}(Y,\cdot)} &{\to}&{\cdots}\\
 & &{\downarrow \ \uparrow} & &{\downarrow\ \uparrow} & &
{\downarrow \  \uparrow} & &\\
{\cdots}&{\to}&{H^i(X,\cdot)}&{\to}&{H^i(X-
f^{-1}(Z),\cdot)} & {\to}& {H^{i+1}_{f^{-1}(Z)}(X,\cdot)}
&{\to}&{\cdots} \end{array} $$
Taking the direct limit of the concerned diagram over
$Z\in Z^1(Y)$ we do get
$$ \begin{array}{ccccccccc}
{\cdots}&{\to}&{H^i(Y,\cdot)}&{\to}&{H^i(K(Y),\cdot)} &
{\to}& {H^{i+1}_{Z^1(Y)}(Y,\cdot)} &{\to}&{\cdots}\\
 & &{\downarrow \ \uparrow} & &{\vcong\ \uparrow} & & {\downarrow \
\uparrow} & &\\
{\cdots}&{\to}&{H^i(X,\cdot)}&{\to}&{H^i(K(X),\cdot)} &
{\to}& {H^{i+1}_{Z^1(X)}(X,\cdot)} &{\to}&{\cdots}
\end{array} $$ Thus, taking the associated sheaves, we
have:  $$ \begin{array}{ccccccccc}
{\cdots}&{\by{zero}}&{\cH^i_{Y}\p{(\cdot)}}&{\to}&{
i_YH^i(K(Y),\p{\cdot})} & {\to}&
{\cH^{i+1}_{Z^{1},Y}\p{(\cdot)}} &{\by{zero}}&{\cdots}\\
 & &{\downarrow\ \uparrow} & &{\vcong\ \uparrow} & &
{\downarrow \  \uparrow}  & &\\
{\cdots}&{\to}&{af_*H^i_{X}\p{(\cdot)}}&{\to}&{af_*H^i(K(X),\cdot)}
& {\to}& {af_*(H^{i+1}_{Z^{1},X}\p{(\cdot)})} &{\to}&{\cdots}
\end{array} $$ and furthermore
$$ \begin{array}{ccccccccc}
{\cdots}&{\to}&{af_*H^i_{X}\p{(\cdot)}}&{\to}&{af_*H^i(K(X),
\cdot)} & {\to}& {af_*(H^{i+1}_{Z^{1},X}\p{(\cdot)})}
&{\to}&{\cdots}\\
& &{\downarrow} & &{\vcong} & & {\downarrow} & &\\
{0}&{\to}&{f_*(\cH^i_{X}\p{(\cdot)}})&{\to}&{
f_*(i_XH^i(K(X),\p{\cdot})} & {\to}&
{f_*(\cH^{i+1}_{Z^{1},X}\p{(\cdot)})} &{\to}&{\cdots}
\end{array} $$
One then obtain by patching the diagram (\ref{diagramma
chiave}). Now because of Lemma~\ref{chiave} and
(\ref{diagramma chiave})  we do get the first claimed
isomorphism $f^{\sharp}:
\cH^*_Y\p{(\cdot)} \by{\simeq} f_*\cH^*_X\p{(\cdot)}$\,.

\subsection{Proving the $K$-theory case}\label{svol2}

We now consider the Quillen $K$-theory of vector
bundles. The proof of Theorem~1 is the analogous of the
previous one by using the Gillet's projection formula and
the Gersten's conjecture.\\
To prove the Theorem~1, arguing as in \S\ref{svol1}, we can
assume $X$ and $Y$ irreducibles, $r=0$ and  $f_{!}:
f_*\cI^{\bul}_{0,X} \by{ \simeq} \cI^{\bul}_{0,Y}$ given
by  $f_{!}(1_X)=1_Y$ (because $K(X)\cong K(Y)$). Hence,
 by (\ref{kpr}), we obtain
 $$\B{array}{ccccc}
& &f_*\cK_{p}(\cO_X) & &\\ &\p{f^{\natural}}\nearrow \ & &\
\searrow \p{\cap_{f_*}1_X}& \\ \cK_{p}(\cO_Y)& & & &
f_*\cI^{\bul}_{p,X}\\ \veq & & & & \ \downarrow \p{f_{!}}
\\ \cK_{p}(\cO_Y) & & \longby{\cap 1_Y} & &
\cI^{\bul}_{p,Y}
\E{array}$$
One defines $f_{\natural}:f_*\cK_{p}(\cO_X)
\to\cK_{p}(\cO_Y)$ in the derived category, as follows:
$$f_{\natural} \df (\ \cap 1_Y)^{-1}\p{\circ} f_{!}
\p{\circ} (\ \cap_{f_*}1_X)$$ Thus:
$f_{\natural}\p{\circ}f^{\natural}=id$ .
(Note: the map $\cap_{f_*}1_X$ is not a quasi-isomorphism
in general, but it induces an isomorphism on homology in
degree zero because $f_*$ is left exact; indeed, taking
$h^0$= the zero homology of a complex, we have the
commutative diagram of sheaves
 $$\B{array}{ccccc}
& &f_*\cK_{p}(\cO_X) & &\\
&\p{f^{\natural}}\nearrow \ & &\ \ \searrow \p{\simeq} & \\
\cK_{p}(\cO_Y)& & & & h^0(f_*\cI^{\bul}_{p,X})\\
\veq & & & & \ \downarrow \p{f_{!}} \\
\cK_{p}(\cO_Y) & & \longby{\simeq} & &
h^0(\cI^{\bul}_{p,Y})
\E{array}$$
and $f_{\natural}\p{\circ}f^{\natural}=id$
between sheaves on $Y$ .)\\
Associated to $f:X\to Y$ proper birational morphism
between regular (irreducible) algebraic schemes, we have
 \B{displaymath}\B{array}{ccc}
\cK_{p}(\cO_Y) & \hookrightarrow & i_YK_p(K(Y)) \\
 \p{f_{\natural}} \uparrow \downarrow \p{f^{\natural}}
 & & \vcong\uparrow \ \\
f_*(\cK_{p}(\cO_X)) & \hookrightarrow &
f_*(i_XK_p(K(X)))
\E{array}\E{displaymath}
so that by the same argument as in \S\ref{H-coda}  we do
get the second claimed isomorphism $f^{\natural}:
\cK_{p}(\cO_Y) \by{\simeq} f_*\cK_{p}(\cO_X)$\,.

\B{rmk}
Assuming the Gersten's conjecture and applying the
above argument one can see that $f^{\natural}$ is an
isomorphism if $f$ is a proper birational morphism between
regular biequidimensional schemes.
\E{rmk}

\section{Homotopy and proto-decomposition}

We maintain the notations and the assumed `cohomology
theory' introduced in the previous Section
(see \S1). Let $\P^n_X$ be the scheme $X\times
_{k}{\rm Proj}\, k [t_0,\ldots,t_n]$; let $\pi_n :
\P^n_X\to X$ denote the canonical projection on $X$ smooth
and equidimensional in $\cat{V}_k$ and assume $k$ perfect.
For any couple of non-negative integers $n\geq m$ let
$j_{(n,m)}$ denote the `Gysin homomorphism'(see
\S\ref{cad1} below)  $$H^{p}(\P^m_X,\cH^{q}\p{(j)})\to
H^{p+n-m}(\P^n_X,\cH^{q+n-m}\p{(j+n-m)})$$
given by the smooth pair $(\P^n_X,\P^m_X)$
of pure codimension $n-m$\,; if $m\geq l$ is another such
couple, i.e. $(\P^m_X,\P^l_X)$ is a pair, then $j_{(n,l)} =
j_{(n,m)}\p{\circ}j_{(m,l)}$ .
Let $\A^1_X$ denote the scheme $X\tensor _{k}k [t]$ and
assume that the cohomology theory satisfies the
following.\\[1pt]

{\bf Homotopy property}. {\em Let $X$ be an algebraic
smooth scheme. The natural morphism $\pi: \A^1_X \to X$
induces an isomorphism $$\pi^*: H^*(X,\cdot)
\by{\simeq} H^*(\A^1_X,\cdot)$$  by pulling-back along
$\pi$.}\\[1pt]

For $\cE$ a locally free sheaf on $X$, $\rank \cE =
n+1$ and $\pi:{\bf V}(\cE) \to X$ the
associated vector bundle, we then get the
isomorphism (see \S\ref{minom} below)
$$H^p_Z(X,\cH^q\p{(j)})\cong H^p_{\pi^{-1}(Z)}({\bf
V}(\cE),\cH^q\p{(j)})$$
pulling back along $\pi$ where $Z\subseteq X$ is
any closed subset. Furthermore, it is now possible to prove
the following Dold-Thom type decomposition.

\B{teor} {\rm (Proto-decomposition)}
Let $X$ be algebraic, equidimensional and smooth over a
perfect field. Assuming the homotopy property above then
there is an isomorphism  $$H^p(\P^n_X,\cH^q_X\p{(j)})\cong
\bigoplus_{i=0}^{n} H^{p-i}(X,\cH^{q-i}_X\p{(j-i)})$$
where every $x\in H^p(\P^n_X,\cH^q_X\p{(j)})$ is written
as $$\pi^*_n(x_{n})+j_{(n,n-1)}\pi^*_{n-1}(x_{n-1})+\cdots
+j_{(n,1)}\pi^*_1(x_{1})+j_{(n,0)}(x_0)$$  for
$x_{n-i}\in H^{p-i}(X,\cH^{q-i}_X\p{(j-i)})$ and
$i=0,\ldots ,n$ .
\E{teor}

\B{rmk} Note that for $\cE$ a locally free sheaf on $X$ ,
$\rank \cE = n+1$, we will obtain the decomposition of
$H^p(\P(\cE),\cH^q_X\p{(j)})$ in Scholium~\ref{Edeco}.
\E{rmk}

Before proving the Theorem~2 we need the following results.

\subsection{Gysin maps for $\cH$-cohomologies}\label{cad1}

The category $\cat{V}^2_k$ is the category of pairs of
algebraic schemes over a perfect field $k$ .
 \B{lemma}{\rm (Purity) } If $(X,Z)$ is a pure smooth pair
in $\cat{V}^2_k$, ${\rm codim}_XZ=c$, then
$H^p_Z(X,\cH^q_X\p{(j)})$ is canonically isomorphic to
$H^{p-c}(Z,\cH^{q-c}_Z\p{(j-c)})$. \label{purity}\E{lemma}
\B{proof} Let $\cR_q^{\bul}\p{(j)}$ denote the arithmetic
resolution of the sheaf $\cH^q\p{(j)}$ on $X$ (resp. on
$Z$) and denote $H^0(X,\cR_q^{\bul}\p{(j)}) \df  {\rm
R}_q^{\bul}(X)\p{(j)}$ (resp. ${\rm R}_q^{\bul}(Z)\p{(j)}$
). Then: \B{slemma}\label{sup}
For $Z \subset X$ of pure codimension c :
$$H^0_Z(X,\cR_q^{\bul}\p{(j)}) \cong {\rm
R}_{q-c}^{\bul}(Z)\p{(j-c)}[-c]$$
\E{slemma}
\B{proof}
Straightforward.
\E{proof}
Since  $\cR_q^{\bul}\p{(j)}$ is a  bounded complex
(graded by codimension) of flasque sheaves the
hypercohomology spectral sequence ($h^n(C^{\bul}) \df$ the
n$^{th}$ homology group of a complex $C^{\bul}$)
$$'E^{r,s}_2 = h^r(H^s_Z(X,\cR_q^{\bul}\p{(j)}))
\Rightarrow \H^{r+s}_Z(X,\cR_q^{\bul}\p{(j)})$$
degenerates to isomorphisms
$$ h^p(H^0_Z(X,\cR_q^{\bul}\p{(j)})) \cong
\H^p_Z(X,\cR_q^{\bul}\p{(j)})$$
Taking account of the Sublemma~\ref{sup}, because of the
(flasque) arithmetic resolutions, we do get a chain of
isomorphisms $$\begin{array}{rcl} H^p_Z(X,\cH^q_X\p{(j)})
& \cong & \H^p_Z(X, \cR_q^{\bul}\p{(j)}) \\ & \cong &
h^p(H^0_Z(X,\cR_q^{\bul}\p{(j)}))\\ &   \cong &
h^p({\rm R}_{q-c}^{\bul -c}(Z)\p{(j-c)}) \\
& \cong & H^{p-c}(Z,\cH^{q-c}_Z\p{(j-c)})
\end{array} $$
The proof of the Lemma~\ref{purity} is complete.
\E{proof}
\B{schol}{\rm (Gysin map) } \label{Gysin}
 Let $(X,Z)\in
\cat{V}^2_k$ be a smooth pair of pure
codimension $c$. There is an homomorphism
$$j_{\pp{(X,Z)}} :H^{p}(Z,\cH^{q}_Z\p{(j)})\to
H^{p+c}(X,\cH^{q+c}_X\p{(j+c)})$$     such that if
$(Z,T)$ is another smooth pair then $j_{\pp{(X,T)}} =
j_{\pp{(X,Z)}}\p{\circ}j_{\pp{(Z,T)}}$ .
\E{schol}
\B{proof}
The map $j_{\pp{(X,Z)}}$ is induced on
cohomology by the composition in the derived category
(by sheafifying the isomorphism in the Sublemma~\ref{sup}
and using the arithmetic resolutions)
$$j_*\cH^{q}_Z\p{(j)}\by{\simeq}j_*\cR_{q,Z}^{\bul}\p{(j)}
\by{\simeq}{\bf \Gamma}_Z\cR_{q+c,X}^{\bul}\p{(j+c)}[c]
\to \cR_{q+c,X}^{\bul}\p{(j+c)}[c] \by{\simeq}
\cH^{q+c}_X\p{(j+c)}[c]$$ where $j:Z\hookrightarrow X$
and ${\bf \Gamma}_Z$ are the sections supported in $Z$ .
The compability simply follows by considering the
resolutions $\cR^{\bul}$ and observing that a global
section of $\cR^{\bul}$ on $T$ can be seen as a section of
$\cR^{\bul}$ on $Z$ supported in $T$, shifted by the
codimension of $T$ in $Z$, etc\ldots , as a section of
$\cR^{\bul}$ on $X$ shifted by ${\rm codim}_TZ + {\rm
codim}_ZX = {\rm codim}_TX$.
 \E{proof}

\subsection{Homotopy for $\cH$-cohomologies}\label{minom}

For $X\in \cat{V}_k$ smooth over $k$ perfect, we recall
(see \cite[6.3]{BO}) that exists a spectral sequence
(`coniveau') $$E^{p,q}_2 = H^p(X,\cH^q\p{(\cdot)})
\Rightarrow H^{p+q}(X,\cdot)$$

\B{lemma}{\rm (Homotopy) } If the functor $H^*(\ ,\cdot)$
has the homotopy property then the functor $H^{\#}_{Zar}(\
,\cH^*\p{(\cdot)})$ has the homotopy property.
\label{homotopy} \E{lemma} \B{proof} Let $\pi : \A^1_X \to X$ be the
given structural flat morphism over $X$ smooth. We will
show that $\pi^*: H^{\#}_{T}(X,\cH^*\p{(\cdot)})\cong
H^{\#}_{\pi^{-1}(T)}(\A^1_X,\cH^*\p{(\cdot)})$ for any
closed subscheme $T\subseteq X$. The proof is divided in
two steps.\\ First step: reducing to the function field
case  $\A^1_{K} \to K$ . This is done using a trick
by Quillen \cite[Prop.4.1]{Q}.
Associated to $\pi$ and $Z\subset T\subset X$ closed
subsets, $U =X-Z$ , we have a map of long exact
sequences:   $$ \begin{array}{ccccccccc}
\cdots & \to & H^p_T(X,\cH^*\p{(\cdot)}) & \to &
H^p_{T\cap U}(U,\cH^*\p{(\cdot)}) & \to &
H^{p+1}_Z(X,\cH^*\p{(\cdot)}) & \to & \cdots \\
 & &{\downarrow} & &{\downarrow} & & {\downarrow} & &\\
{\cdots}&{\to}&
H^p_{\pi^{-1}(T)}(\A^1_X,\cH^*\p{(\cdot)}) &{\to}&
H^p_{\pi^{-1}(T\cap U)}(\A^1_U,\cH^*\p{(\cdot)}) & {\to}&
H^{p+1}_{\pi^{-1}(Z)}(\A^1_X,\cH^*\p{(\cdot)})
&{\to}&{\cdots} \end{array} $$
By the five lemma the induced (middle vertical)
homomorphisms $$H^p_T(X,\cH^*\p{(\cdot)}) \to
H^p_{\pi^{-1}(T)}(\A^1_X,\cH^*\p{(\cdot)})$$ are isomorphisms (all $p$)
if the others vertical maps are. Using noetherian
induction we can assume $H^p_Z(X,\cH^*\p{(\cdot)}) \to
H^p_{\pi^{-1}(Z)}(\A^1_X,\cH^*\p{(\cdot)})$ to be an
isomorphism for all closed subsets $Z\neq T$ and all
$p\geq 0$. We can also suppose $X$ irreducible. Taking the
direct limit over all closed proper subschemes $Z$ of $T$
we can also assume that $T$ is integral of codimension
$t$. Thus by local purity we are left to show that
\B{eqnarray}\label{field} \limdir{U = X-Z} H^p_{T\cap
U}(U,\cH^*\p{(\cdot)}) & \cong & H^{p-t}(\Sp
(K(T)),\cH^{*-t}\p{(\cdot-t)})  \nonumber \\ & &\ \ \ \ \ \
\downarrow  \\ \limdir{U = X-Z} H^p_{\pi^{-1}(T\cap
U)}(\A^1_U,\cH^*\p{(\cdot)}) & \cong &
H^{p-t}(\A^1_{K(T)},\cH^{*-t}\p{(\cdot-t)}) \nonumber
\E{eqnarray} is an isomorphism for all $p$. (Note: the
horizontal isomorphisms in (\ref{field}) are obtained by
continuity of the arithmetic resolution of the sheaf
$\cH^*\p{(\cdot)}$ ).\\  Second step: proving the function
field case $K=K(X)$. Having defined $$H^*(\A^1_{K},\cdot)
\df \limdir{U\subset X}  H^*(\A^1_{U},\cdot)$$ by
continuity of the coniveau spectral sequence we do have
$$E^{p,q}_2 = H^p(\A^1_{K},\cH^q\p{(\cdot)}) \Rightarrow
H^{p+q}(\A^1_{K},\cdot)$$ and $E^{p,q}_2 =0$ if $p>1={\rm
dim}\, \A^1_{K}$, thus all the differentials are zero
which just yields short exact sequences     $$0 \to
H^1(\A^1_{K},\cH^{q-1}\p{(\cdot)})\to
H^q(\A^1_{K},\cdot)\to  H^0(\A^1_{K},\cH^q\p{(\cdot)}) \to
0 $$ Associated to the flat map $\A^1_K \to K$ we have  a
commutative diagram
\B{equation}\begin{array}{ccccccccc}\label{split} 0 & \to
& H^1(\A^1_{K},\cH^{q-1}\p{(\cdot)}) & \to &
H^q(\A^1_{K},\cdot) & \to & H^0(\A^1_{K},\cH^q\p{(\cdot)})
 & \to & 0\\
 & & & &{\cong\uparrow\;\;} & & {\downarrow \ \uparrow}
& &\\
& & & & H^q(K,\cdot) & \by{\simeq} &
H^0(K,\cH^q\p{(\cdot)}) & & \end{array}
\E{equation}
where: because $H^*(\ ,\cdot)$ has the homotopy property
then the middle vertical map is an isomorphism while the
arrow $H^0(\A^1_{K},\cH^q\p{(\cdot)})\to H^q(K,\cdot)$ is
the evaluation at any $K$-rational point (cf.
\cite[Proof of 2.5]{BR}). From (\ref{split}) it follows
that $H^1(\A^1_{K},\cH^{q-1}\p{(\cdot)}) = 0$ hence the
required isomorphism in (\ref{field}) is given by:
$ H^0(K,\cH^q\p{(\cdot)}) \cong
H^0(\A^1_{K},\cH^q\p{(\cdot)})$
\E{proof}

\B{schol}\label{hind}
Let $\A^n_X$ denote the n$^{th}$ affine space over
$X$ smooth (i.e. the scheme $X\tensor
_{k}k [t_1,\ldots,t_n]$ ) and let $\pi : \A^n_X \to X$ be
the natural projection. Then $\pi$ induces an isomorphism
$$H^{\#}_{\pi^{-1}(Z)}(\A^n_X,\cH^*\p{(\cdot)})\cong
 H^{\#}_Z(X,\cH^*\p{(\cdot)})$$  \E{schol}
\B{proof}
By induction from the Lemma~\ref{homotopy}.
\E{proof}

\B{cor}\label{Ehomo} Let  $\cE$ be a locally free sheaf on
$X$ smooth and $\pi:{\bf V}(\cE) \to X$ the
associated vector bundle, we then get the
isomorphism $$\pi^*: H^p_Z(X,\cH^q\p{(j)})\cong
H^p_{\pi^{-1}(Z)}({\bf V}(\cE),\cH^q\p{(j)})$$
\E{cor}
\B{proof} By reduction to open Zariski neighborhoods on
which $\cE$ is free and noetherian induction (cf.
the proof of Lemma~\ref{homotopy}).
\E{proof}

\subsection{Proof of the proto-decomposition
Theorem}\label{svol3}

The proof of Theorem~2 is by induction on $n$. For $n=0$,
$\P^0_X \cong X$, hence the induction starts: one consider
an `hyperplane at infinity' $\infty$ in $\P_X^n$ so that
$\infty \cong \P_X^{n-1}$ and $\P_X^n -
\infty\cong\A_X^n$. There is a standard long exact
sequence of Zariski cohomology groups
$$H^{p-1}(\A_X^n,\cH^q\p{(j)}) \to H^{p}_{\infty}(\P_X^n,
\cH^q\p{(j)}) \to H^{p}(\P_X^n,\cH^q\p{(j)}) \to
H^{p}(\A_X^n,\cH^q\p{(j)})$$ Let $\pi_n :\P_X^n\to X$ be
the projection. By homotopy (see \S\ref{minom}) the
restriction $\pi_n \mid_{\pp{\A_X^n}} : \A_X^n \to X$
induces a splitting of the previous long exact sequence,
given by the commutative square $$\begin{array}{ccc}
H^{p}(\P_X^n,\cH^q\p{(j)})&\to
&H^{p}(\A_X^n,\cH^q\p{(j)})\\
\p{\pi_n^*}\uparrow & &\uparrow \cong \ \\
 H^{p}( X,\cH^q\p{(j)})& =  &H^{p}( X,\cH^q\p{(j)})
\end{array}
$$
By the purity Lemma~\ref{purity} we have an isomorphism
$$H^{p}_{\infty}(\P_X^n,\cH^q\p{(j)})\cong
H^{p-1}(\P_X^{n-1},\cH^{q-1}\p{(j-1)})$$
and, by composing, the Gysin map (see \S\ref{cad1})
$$j_{(n,n-1)}:H^{p-1}(\P^{n-1}_X,\cH^{q-1}\p{(j-1)})\to
H^{p}(\P^n_X,\cH^{q}\p{(j)})$$
So one can split the long exact sequence into short exact
sequences $$0 \to H^{p-1}(\P^{n-1}_X,\cH^{q-1}\p{(j-1)})
\by{j_{(n,n-1)}} H^{p}(\P^n_X,\cH^{q}\p{(j)})
\stackrel{\stackrel{\pi_n^*}{\leftarrow}}{\to}
H^{p}(X,\cH^q\p{(j)}) \to 0 $$
Thus we do get the formula: $$H^p(\P
_X^n,\cH^q\p{(j)})\cong H^{p}( X,\cH^q\p{(j)}) \oplus
H^{p-1}(\P_X^{n-1},\cH^{q-1}\p{(j-1)})$$ which does the
induction's step: an element $x\in
H^p(\P^n_X,\cH^q_X\p{(j)})$ is written as
$\pi^*_n(x_n)+j_{(n,n-1)}(x')$  for  $x_n\in
H^{p}(X,\cH^{q}_X\p{(j)})$ and
$x'\in H^{p-1}(\P_X^{n-1},\cH^{q-1}\p{(j-1)})$; this last
$x'$ because of the inductive hypothesis is written as
$$\pi^*_{n-1}(x_{n-1})+j_{(n-1,n-2)}\pi^*_{n-2}(x_{n-2})+
\cdots +j_{(n-1,1)}\pi^*_1(x_{1})+j_{(n-1,0)}(x_0)$$
for  $x_{n-1-i}\in H^{p-1-i}(X,\cH^{q-1-i}_X\p{(j-1-i)})$
and $i=0,\ldots ,n-1$ and because of the compatibility of
the Gysin homomorphisms (see Scholium~\ref{Gysin}) applying
$j_{(n,n-1)}$ we are done.

\section{Cap products}

This section is devoted to construct a cap product
between algebraic cycles and $\cH$-cohomology classes.

\subsection{Sophisticated Poincar\'e duality theories}

Let assume we are given a cohomology theory
$H^*(\p{\cdot})$ and a homology theory $H_*(\p{\cdot})$ on
$\cat{V}_k$ satisfying the Bloch-Ogus axioms
\cite[1.1-1.2]{BO}. Furthermore, we let assume the
existence of a {\it sophisticated}\, cap-product with
supports i.e. for all $(X,Z),(X,Y)\in\cat{V}^2_k$ a
pairing   $$\p{\cap_{Y,Z}}: H_n(Y,\p{m})\tensor
H_Z^q(X,\p{s})\to H_{n-q}(Y\cap Z, \p{m-s})$$ which
satisfies the following axioms: \B{description} \item[A1]
$\p{\cap}$ is natural with respect to \'etale maps (or
just open Zariski immersions according with
\cite[1.4.2]{BO}) of pairs in $\cat{V}^2_k$.  \item[A2] If
$(X,T),(T,Z)$ are pairs in $\cat{V}^2_k$ then the
following diagram
 \B{displaymath}\B{array}{ccc}
H_n(Z,\p{m})\tensor H^q(X,\p{s}) & \longby{\cap} &
H_{n-q}(Z,\p{m-s})\\
 \downarrow & &  \downarrow \\
H_n(T,\p{m})\tensor H^q(X,\p{s})  & \longby{\cap} &
H_{n-q}(T,\p{m-s})
\E{array}\E{displaymath}
commutes.
\item[A3] For $(X,T),(T,Z)$ pairs in
$\cat{V}^2_k$ let $U=X-Z$ and let $j: U \to X$ be the
inclusion. Let denote $\partial:
H_n(T\cap U,\cdot)\to H_{n-1}(Z,\cdot)$ the
boundary map in the long exact sequence (\ref{hloc}) of
homology groups. Then the following diagram
\B{displaymath}\B{array}{c}
\ \hspace{20pt} H_n(T\cap U,\cdot)\tensor
H^q(U,\cdot) \\ \p{id\tensor j^*}\nearrow \hspace{50pt}
\searrow\p{\cap}\\H_n(T\cap U,\cdot)\tensor
H^q(X,\cdot)\hspace{40pt} H_{n-q}(T\cap U,\cdot) \\
\p{\partial\tensor id}\downarrow \hspace{97pt} \downarrow
\p{\partial} \\ H_{n-1}(Z,\cdot)\tensor H^q(X,\cdot)
\hspace{10pt} \longby{\cap} \hspace{10pt}
H_{n-q-1}(Z,\cdot)\E{array} \E{displaymath}
 commutes i.e. we have the equation:
\B{equation} \partial
(y\p{\cap}j^*(x))=\partial (y)\p{\cap}x
\E{equation}
for $y\in H_n(T\cap U,\cdot)$ and $x\in
H^q(X,\cdot)$.
\item[A4 {\it (Projection Formula)}] Let $f:X'\to X$ be a
proper morphism in $\cat{V}_k$. For $(X,Y)$ and $(X,Z)$ let
$Y'=f^{-1}(Y)$ and $Z'=f^{-1}(Z)$. The following diagram

$$\B{array}{c}
\hspace{20pt} H_{n}(Y',\p{m})\tensor
H^{q}_{Z'}(X',\p{s})\\ \p{id\tensor
f^{*}}\nearrow \ \hspace{40pt} \searrow
\p{\cap}  \\H_{n}(Y',\p{m})\tensor
H^{q}_Z(X,\p{s})\hspace{60pt}
H_{n-q}(Y'\cap Z',\p{m-s})\\  \p{f_{*}\tensor
id}\downarrow \hspace{121pt}\downarrow \p{f_{*}} \\
H_{n}(Y,\p{m})\tensor
H^{q}_Z(X,\p{s})\hspace{10pt} \longby{\cap}
\hspace{10pt}H_{n-q}(Y\cap Z,\p{m-s})
\E{array}$$
commutes.
\E{description}

By the way, for $(H^*,H_*)$ as above, we have the
following `projection formula':

\B{schol}\label{A4} Let $f:Y\to X$ be a proper morphism
in $\cat{V}_k$.  Let $T$ be any closed subset of $Y$ and
let $f(T)=Z$. Then the following diagram

\B{equation}\B{array}{c}
\ \hspace{20pt}
H_{n}(T,\p{m})\tensor H^q(Y,\p{s})
\\ \p{id\tensor f^*}\nearrow \hspace{50pt}
\searrow\p{\cap}\\  H_{n}(T,\p{m})\tensor
H^q(X,\p{s})\hspace{40pt} H_{n-q}(T,\p{m-s}) \\
\p{f_{*}\tensor id}\downarrow \hspace{105pt} \downarrow
\p{f_{*}} \\H_{n}(Z,\p{m})\tensor H^q(X,\p{s})
\hspace{10pt} \longby{\cap} \hspace{10pt}
H_{n-q}(Z,\p{m-s}) \E{array}  \E{equation}
commutes.
\E{schol}
\B{proof} This is a simple consequence of A4 by observing
that $T\into f^{-1}(Z)$.
\E{proof}

\B{defi} We will say that $(H^*,H_*)$ is a {\it
sophisticated} Poincar\'e duality theory with supports if
the axioms A1--A4 are satisfied and Poincar\'e duality
holds i.e. the Bloch-Ogus axioms \cite[1.3.4-5 and
7.1.2]{BO} are satisfied (see \S\ref{inter1}).
\E{defi}

\subsection{$\cH$-cap product}\label{H-cap}

Associated with the homology theory $H_*$, for
$X\in\cat{V}_k$ possibly singular, we have a niveau
spectral sequence (cf. \cite[Prop.3.7]{BO}) $$E_{a,b}^1 =
\coprod_{x\in X_a}^{} H_{a+b}(x,\p{\cdot}) \Rightarrow
H_{a+b}(X,\cdot)$$ which is covariant for proper morphisms
and contravariant for \'etale maps. Let denote  ${\rm
Q}^n_{\bul}(X)\p{(m)}$ the (homological) complex
$E_{\bul,n}^1\p{(m)}$.

\B{prop} \label{hpairing} Let $H^*$ and $H_*$ be
cohomological and homological functors satisfying the
axioms {\rm A1--A3} above. For $X\in\cat{V}_k$ there is a
pairing of complexes $${\rm Q}^n_{\bul}(X)\p{(m)}\tensor
H^q(X,\p{s})\to {\rm Q}^{n-q}_{\bul}(X)\p{(m-s)}$$
contravariant w.r.t. \'etale maps.  \E{prop}
\B{proof}
Let $Z\subset T\subseteq X$ be closed subsets of $X$,
dim$T\leq a$, dim$Z\leq a-1$ and let $U=X-Z$; thus by
restriction to $U$ and cap-product we do have a pairing
associated to such pairs $Z\subseteq T$:  $$H_i(T\cap
U,\p{j})\tensor H^q(X,\p{s}) \to H_{i-q}(T\cap U,\p{j-s})
$$  i.e. $t\p{\tensor}x \leadsto t\p{\cap}j^*(x)$ where $j:
U\hookrightarrow X$.
By taking the direct limit over
such pairs (this makes sense because of A1--A2)  we do
have a pairing  $$\coprod_{x\in X_a}^{}
H_{n+a}(x,\p{m})\tensor H^q(X,\p{s}) \to
\coprod_{x\in X_a}^{} H_{n-q+a}(x,\p{m-s})$$
We need to check compatibility with the differentials of
${\rm Q}^n_{\bul}(X)\p{(m)}$. Because of A2 we have a
pairing  $H_{i}(Z_{a},\p{j})\tensor
H^q(X,\p{s})\to H_{i-q}(Z_{a},\p{j-s})$ where
$H_{i}(Z_{a},\p{j})  \df \limdir{\pp{T\subset X}}
H_i(T,\p{j})$, $T$ as above. Because of A3 and
limit arguments the following diagram
\B{displaymath} \begin{array}{ccc}
{\displaystyle\coprod_{x\in X_a}^{}
H_{i}(x,\p{j})\tensor H^q(X,\p{s})} &
\longby{\partial\tensor id} &
H_{i-1}(Z_{a-1},\p{j}) \tensor H^q(X,\p{s})\\
{\downarrow \p{\cap}} & &{\downarrow \p{\cap} }\\
{\displaystyle\coprod_{x\in X_a}^{}
H_{i-q}(x,\p{j-s})} &\to& H_{i-q-1}(Z_{a-1},\p{j-s})
 \end{array}
 \E{displaymath}
commutes (indeed $ \partial (t\p{\cap}j^*(x))=\partial
(t)\p{\cap}x$ ) and A1 implies that the following
 \B{displaymath} \begin{array}{ccc}
H_{i-1}(Z_{a-1},\p{j}) \tensor H^q(X,\p{s}) &\to &
{\displaystyle\coprod_{x\in X_{a-1}}^{}
H_{i-1}(x,\p{j}) \tensor H^q(X,\p{s}) }\\
{\downarrow \p{\cap}} & &{\downarrow \p{\cap}}\\
H_{i-q-1}(Z_{a-1},\p{j-s})&\to &{\displaystyle
\coprod_{x\in X_{a-1}}^{}
H_{i-q-1}(x,\p{j-s})}
 \end{array}
 \E{displaymath}
commutes. By construction, the differential is the
composition of  $${\rm Q}^n_{a}(X)\p{(m)}\to
H_{n+a-1}(Z_{a-1},\p{m})\to{\rm Q}^n_{a-1}(X)\p{(m)}.$$
Thus the result.
\E{proof}

\B{defi} For $X\in\cat{V}_k$, by taking
associated sheaves for the Zariski topology of the pairing
above, we get a pairing
$$\p{\cap}_{\cH}:{\cQ}^n_{\bul ,X}\p{(m)}\tensor
\cH^q_X\p{(s)}\to {\cQ}^{n-q}_{\bul ,X}\p{(m-s)}$$ which
we call {\it $\cH$-cap-product} on $X$.
\E{defi}

\subsection{Projection formula}

Let $f:Y\to X$ be a proper morphism in $\cat{V}_k$.
We do have a map of niveau
spectral sequences  $$E_{p,q}^1\p{(r)}(Y) \to
E_{p,q}^1\p{(r)}(X)$$
which takes $y\in Y_{p}$ to $f(y)$ if dim
$\overline{\{f(y)\}} =$ dim $\overline{\{y\}}$ zero
otherwise and maps $H_i(y)$ to $H_i(f(y))$. Thus by
sheafifying it for the Zariski topology we obtain a map
$$f_{\sharp}: f_*{\cQ}^n_{\bul ,Y}\p{(m)} \to {\cQ}^n_{\bul
,X}\p{(m)}$$
of complexes of sheaves on $X$.
 \B{prop} The following diagram:
\B{equation}
\B{array}{c}\label{hpr}
\hspace{10pt} f_*{\cQ}^n_{\bul ,Y}\p{(m)}\tensor
f_*\cH^{q}\p{(s)}\\ \p{id\tensor f^{\sharp}}\nearrow \
\hspace{40pt} \searrow \p{f_*\cap_{{\cH}}}  \\
f_*{\cQ}^n_{\bul ,Y}\p{(m)}\tensor
\cH^{q}\p{(s)}\hspace{60pt}
f_*{\cQ}^{n-p}_{\bul ,Y}\p{(m-s)}\\  \p{f_{\sharp}\tensor
id}\downarrow \hspace{105pt}\downarrow \p{f_{\sharp}} \\
{\cQ}^n_{\bul ,X}\p{(m)}\tensor
\cH^{q}\p{(s)} \hspace{10pt} \longby{\cap_{\cH}}
\hspace{10pt}{\cQ}^{n-q}_{\bul ,X}\p{(m-s)}
\E{array}
\E{equation}
commutes.
\E{prop}

\B{proof} The commutative diagram above will be otbained
from the following:
$$
\B{array}{c}
\hspace{20pt} {\rm Q}^n_{\bul}(f^{-1}(U))\p{(m)}\tensor
H^{q}(f^{-1}(U),\p{s})\\ \p{id\tensor f^{*}}\nearrow \
\hspace{40pt} \searrow \p{\cap}  \\
{\rm Q}^n_{\bul}(f^{-1}(U))\p{(m)}\tensor
H^{q}(U,\p{s})\hspace{60pt}
{\rm Q}^{n-q}_{\bul}(f^{-1}(U))\p{(m-s)}\\  \p{f_{*}
\tensor id}\downarrow \hspace{105pt}\downarrow
\p{f_{*}} \\
{\rm Q}^{n}_{\bul}(U)\p{(m)}\tensor
H^{q}(U,\p{s}) \hspace{10pt} \longby{\cap}  \hspace{10pt}
{\rm Q}^{n-q}_{\bul}(U)\p{(m-s)}
\E{array} $$
where $U\subset X$ is any Zariski open subset of $X$, by
taking associated sheaves on $X_{Zar}$.\\
Moreover it sufficies to prove the case of $U=X$.\\
Let $\overline{\{y\}}\subset Y$ such that
$y\in Y_{p}$ and $f(y)\in X_{p}$.
The maps
$$f_{!,y}: \limdir{\pp{V\subset
Y}}H_*({\overline{\{y\}}\cap V}) \to \limdir{\pp{U\subset
X}}H_*({\overline{\{f(y)\}}\cap U})$$
are defined by mapping the elements of
$H_*({\overline{\{y\}}\cap (Y-T)})$ by restriction and
the induced maps
$$f_{!,y}: H_*({\overline{\{y\}}\cap (Y-f^{-1}(Z))}) \to
H_*({\overline{\{f(y)\}}\cap (X-Z)})$$
where $Z=f(T)$ (Note: compabilities are ensured by
\cite[1.2.2]{BO}).  By the definition of the pairing in
Proposition \ref{hpairing} we are left to show that the
following diagram

\B{equation}\B{array}{c} \label{hopen}
\ \hspace{20pt}
H_{*}(\overline{\{y\}}\cap V,\cdot)\tensor H^q(V,\p{s})
\\ \p{id\tensor f^*}\nearrow \hspace{50pt}
\searrow\p{\cap}\\ H_{*}(\overline{\{y\}}\cap V,\cdot)
\tensor H^q(U,\p{s})\hspace{40pt}
H_{*-q}(\overline{\{y\}}\cap V,\cdot \p{-s}) \\
\p{f_{!,y}\tensor id}\downarrow \hspace{105pt}
\downarrow  \p{f_{!,y}} \\ H_{*}(\overline{\{f(y)\}}\cap
U,\cdot)\tensor H^q(U,\cdot) \hspace{10pt} \longby{\cap}
\hspace{10pt} H_{*-q}(\overline{\{f(y)\}}\cap
U,\cdot \p{-s})
\E{array}  \E{equation}
commutes where $f^{-1}(U) = V$. Since
$f(\overline{\{y\}}\cap f^{-1}(U))=
\overline{\{f(y)\}}\cap U$ the diagram (\ref{hopen})
commutes because of the projection formula in the
Scholium~\ref{A4}.
\E{proof}

\B{lemma} If $i: Z\into X$ is a closed embedding then the
canonical map $i_{\sharp}: i_*{\cQ}^n_{\bul}\p{(m)} \to
{\cQ}^n_{\bul}\p{(m)}$ admits a factorisation by a
quasi-isomorphism $i_Z: i_*{\cQ}^n_{\bul}\p{(m)}
\by{\sim}{\bf \Gamma}_Z{\cQ}^n_{\bul}\p{(m)}$ and the
natural map ${\bf \Gamma}_Z{\cQ}^n_{\bul}\p{(m)} \into
{\cQ}^n_{\bul}\p{(m)}$
\E{lemma}
\B{proof} Clear (cf. the proof of the Scholium~3.3).
\E{proof}

\B{cor}
For $i: Z\into X$ as above the following diagram

\B{equation}
\B{array}{c}\label{hipr}
\hspace{10pt} i_*{\cQ}^n_{\bul}\p{(m)} \tensor
i_*\cH^{q}\p{(s)}\\ \p{id\tensor i^{\sharp}}\nearrow \
\hspace{40pt} \searrow \p{i_*\cap_{{\cH}}}  \\
i_*{\cQ}^n_{\bul}\p{(m)} \tensor
\cH^{q}\p{(s)}\hspace{60pt}
i_*{\cQ}^{n-q}_{\bul}\p{(m-s)} \\  \p{i_{Z}\tensor
id}\downarrow\wr \hspace{105pt}\wr\downarrow \p{i_{Z}} \\
{\bf \Gamma}_Z{\cQ}^n_{\bul}\p{(m)} \tensor
\cH^{q}\p{(s)} \hspace{10pt} \longby{\cap_{\cH}}
\hspace{10pt}{\bf \Gamma}_Z{\cQ}^{n-q}_{\bul}\p{(m-s)}
\E{array}
\E{equation}
commutes.
 \E{cor}
\B{proof} This follows by the factorisation of
$i_{\sharp}$ and the projection formula (cf. the Lemma
and the Proposition above).
 \E{proof}

\subsection{Algebraic cycles}\label{algcyc}

We are now going to consider the cycle group naturally
involved with a fixed theory $(H^*,H_*)$ on $\cat{V}_k$.
To this aim we need to assume a `dimension axiom' (cf.
\cite[7.1.1]{BO}). Let assume that our cohomology
theory $H^*$ takes values in a fixed category of
$\La$-modules where $\La = H^0(k,\p{0})$ is a commutative
ring with 1.\\
\B{defi} We will say that $(H^*,H_*)$ satisfies the {\it
dimension axiom}\, when the following properties
\B{description}   \item[A5] If dim$X\leq d$ then
$H_i(X,\p{m})=0$ for $i>2d$. \item[A6] If $X$ is
irreducible then the canonical map $\lambda^*: \La
\by{\sim} H^0(X,\p{0})$, induced by the structural
morphism $\lambda: X\to k$, is an isomorphism.
 \E{description}
hold for any $X\in \cat{V}_k$. We will say that
$(H^*,H_*)$ satisfies the {\it point axiom}\, if the
properties A5-A6 above just hold locally, at the generic
point of any integral subvariety of each $X\in \cat{V}_k$.
\E{defi}
\B{rmk} Clearly the dimension axiom implies the point
axiom.\E{rmk}
Thus: if $X$ is reduced and $\Sigma$ is its singular
locus then $$H_{i}(X,\p{m}) \cong H_{i}(X-\Sigma ,\p{m})$$
for $i\geq 2d$ by applying A5 to the long exact sequence
of homology groups. In particular, if $X$ is integral of
dimension $d$ then by A6 $$H_{2d}(X,\p{d}) \cong
H_{2d}(X-\Sigma ,\p{d}) \cong H^0(X-\Sigma ,\p{0}) \cong
\La $$
Regarding ${\cQ}^n_{\bul}\p{(m)}$ as a (cohomological)
complex of flasque sheaves graded by negative degrees we
do have
$$ \H^{-p}(X,{\cQ}^n_{\bul}\p{(m)})\cong \frac{
{\rm ker} (\coprod_{x\in X_p}^{} H_{n+p}(x,\p{m})  \to
\coprod_{x\in X_{p-1}}^{} H_{n+p-1}(x,\p{m}) )}{{\rm im}
(\coprod_{x\in X_{p+1}}^{} H_{n+p+1}(x,\p{m}) \to
\coprod_{x\in X_p}^{} H_{n+p}(x,\p{m}))}$$
In particular, for $n=p=m$  and the `dimension axiom'
above we do have
$$\H^{-n}(X,{\cQ}^n_{\bul}\p{(n)})\cong {\rm coker}
(\coprod_{x\in X_{n+1}}^{} H_{2n+1}(x,\p{n}) \to
\coprod_{x\in X_n}^{} \La)$$
where ${\displaystyle\coprod_{x\in X_n}^{} \La}$ is the
$\La$-module of algebraic cycles of dimension $n$ in
$X$.\\

\B{defi} Let $(H^*,H_*)$ be a theory satisfying the
point axiom. We will denote
$$C_n(X;\La)\df \H^{-n}(X,{\cQ}^n_{\bul}\p{(n)})$$
the corresponding group of $n$-dimensional algebraic
$\La$-cycles modulo the equivalence relation given by
$$ {\rm im}(\coprod_{x\in X_{n+1}}^{} H_{2n+1}(x,\p{n})\to
\coprod_{x\in X_n}^{} \La)$$
the image of the differential of the niveau spectral
sequence.
 \E{defi}
\hfill\\

For $i: Z\into X$ a closed embedding we clearly do have an
isomorphism
$$C_n(Z;\La) \cong \H^{-n}_Z(X,{\cQ}^n_{\bul}\p{(n)}).$$
Thus, by taking hypercohomology with supports, the
$\cH$-cap-product yields a cap product
\B{equation}
C_n(Z;\La) \tensor H^p_Y(X,\cH^p\p{(p)})\to
C_{n-p}(Z\cap Y;\La)
\E{equation}
Because of the projection formula (\ref{hpr}) this cap
product will have a projection formula  as well.

\section{Cup products}

We are now going to show that the cup-product in
cohomology give us a nice intersection theory for
$\cH$-cohomology. To this aim we will assume the algebraic
schemes in $\cat{V}_k$ to be equidimensionals and $k$ to
be a perfect field. The main results hold true just
for non-singular varieties nevertheless we will not
assume this hypothesis a priori.

\subsection{Multiplicative Poincar\'e duality theories}

Let assume we are given a twisted cohomology theory
$H^*(\p{\cdot})$ on $\cat{V}_k$ in the sense of
Bloch-Ogus \cite[1.1]{BO}. Furthermore, we want to assume
the existence of a cup-product i.e. for all
$(X,Z),(X,Y)\in\cat{V}^2_k$ an associative anticommutative
pairing   $$\p{\cup_{Y,Z}}: H_Y^p(X,\p{r})\tensor
H_Z^q(X,\p{s})\to H^{p+q}_{Y\cap Z}(X,\p{r+s})$$ which
satisfies the following axioms: \B{description}
\item[$\forall$1] $\p{\cup}$ is natural with respect to
pairs in $\cat{V}^2_k$
\item[$\forall$2] If $(X,T),(T,Z)$ are pairs in $\cat{V}^2_k$
then the following diagram
 \B{displaymath}\B{array}{ccc}
H_Z^p(X,\p{r})\tensor H^q(X,\p{s}) & \longby{\cup} &
H^{p+q}_{Z}(X,\p{r+s})\\
 \downarrow & &  \downarrow \\
H_T^p(X,\p{r})\tensor H^q(X,\p{s})  & \longby{\cup} &
H^{p+q}_{T}(X,\p{r+s})
\E{array}\E{displaymath}
commutes.
\item[$\forall$3] For $(X,T),(T,Z)$ pairs in
$\cat{V}^2_k$ let $U=X-Z$ and let $j: U \to X$ be the
inclusion. Let denote $\partial:
H_{T\cap U}^p(U,\cdot)\to H_Z^{p+1}(X,\cdot)$ the
boundary map in the long exact sequence (\ref{loc}) of
cohomology with supports. Then the following diagram
\B{displaymath}\B{array}{c}
\ \hspace{20pt} H_{T\cap U}^p(U,\cdot)\tensor
H^q(U,\cdot) \\ \p{id\tensor j^*}\nearrow \hspace{50pt}
\searrow\p{\cup}\\ H_{T\cap U}^{p}(U,\cdot)\tensor
H^q(X,\cdot)\hspace{40pt} H_{T\cap U}^{p+q}(U,\cdot) \\
\p{\partial\tensor id}\downarrow \hspace{97pt} \downarrow
\p{\partial} \\ H_Z^{p+1}(X,\cdot)\tensor H^q(X,\cdot)
\hspace{10pt} \longby{\cup} \hspace{10pt}
H_Z^{p+q+1}(X,\cdot) \E{array} \E{displaymath}
 commutes i.e. we have the equation:
\B{equation} \partial
(y\p{\cup}j^*(x))=\partial (y)\p{\cup}x
\E{equation}
for $y\in H_{T\cap U}^{p}(U,\cdot)$ and $x\in
H^q(X,\cdot)$.
 \E{description}

\B{defi} A twisted cohomology theory with supports
$H^*$ has a {\it cup-product}\, if there is a pairing
$\p{\cup_{Y,Z}}:
H_Y^p(X,\p{r})\tensor H_Z^q(X,\p{s})\to
H^{p+q}_{Y\cap Z}(X,\p{r+s})$ which satisfies the axioms
$\forall$1--$\forall$3 listed above.
\E{defi}
\hfill\\[4pt]

If furthermore $(H^*,H_*)$ is a Poincar\'e duality theory
we let assume that the  following compatibility between
cap and cup products holds:
 \B{description}
\item[$\forall$4] For $X$ smooth of dimension $d$ let $\eta_X\in
 H_{2d}(X,\p{d})$ be the fundamental class. Then the
following diagram, where $q+j=2d, s+n=d$,
\B{displaymath}\B{array}{ccc}
 H^q(X,\p{s})\tensor H_Z^p(X,\p{r})  & \longby{\cup} &
H^{p+q}_{Z}(X,\p{r+s})\\  \p{\eta_X\cap -\tensor id}
\downarrow & &  \downarrow\p{\eta_X\cap -} \\H_j(X,\p{n})
\tensor H_Z^p(X,\p{r}) & \longby{\cap} &
H_{j-p}(Z,\p{n-r})  \E{array}\E{displaymath}
commutes i.e. we have the equation:
\B{equation}\label{cap=cup}
(\eta_X\p{\cap}x)\p{\cap}z=
\eta_X\p{\cap}(x\p{\cup}z)
\E{equation}
for $x\in H^q(X,\p{s})$ and $z\in  H_Z^p(X,\p{r})$.
\E{description}
\hfill\\[4pt]

Let $f:Y\to X$ be a proper map of {\it smooth}\,
equidimensional algebraic schemes. Let dim $X$ = $\delta$
and  dim $Y$ = $d$. Let $\rho =\delta - d$. For
$Z\subseteq X$ a closed subset there are maps  $$f_{!}:
H_{2d-i}(f^{-1}(Z),\p{d-i})\to
H_{2d-i}(Z,\p{d-i})$$
Because of Poincar\'e duality $f_{!}$ induces a Gysin map
$$f_{*}:H^{i}_{f^{-1}(Z)}(Y,\p{i})\to
H^{i+2\rho}_{Z}(X,\p{i+\rho})$$
which is uniquely determined by the equation
\B{equation}\label{cov}
f_!(\eta_Y\p{\cap}y)=\eta_X\p{\cap}f_*(y)
\E{equation}
for $y\in H^{i}_{f^{-1}(Z)}(Y,\p{i})$. Thus $H^*$ is a
covariant functor w.r.t. proper maps of pairs $(X,Z)$ where
$Z$ is a closed subset of $X$ smooth: indeed $(1_X)_*=
id$ because of  $(1_X)_!=id$ and $(f\p{\circ}g)_* =
f_*\p{\circ} g_*$ because of
$(f\p{\circ}g)_!(\eta\p{\cap}-)=f_!(g_!(\eta\p{\cap}-))=
f_!(\eta\p{\cap}g_*(-))=\eta\p{\cap}f_*(g_*(-))$.\\
The projection formula by \cite{BO} w.r.t. the cap product
give us, via $\forall$4, the following projection formula:

$$\B{array}{c}
\hspace{20pt} H^{p}(Y,\p{r})\tensor
H^{q}_{f^{-1}(Z)}(Y,\p{s})\\ \p{id\tensor
f^{*}}\nearrow \ \hspace{40pt} \searrow
\p{\cup}  \\ H^{p}(Y,\p{r})\tensor
H^{q}_Z(X,\p{s})\hspace{60pt}
H^{p+q}_{f^{-1}(Z)}(Y,\p{r+s})\\  \p{f_{*}\tensor
id}\downarrow \hspace{121pt}\downarrow \p{f_{*}} \\
H^{p+2\rho}(X,\p{r+\rho})\tensor
H^{q}_Z(X,\p{s})\hspace{10pt} \longby{\cup}
\hspace{10pt} H^{p+q+2\rho}_{Z}(X,\p{r+s+\rho})
\E{array}$$
Indeed we have:
\B{center}
\parbox{3in}{$f_!(\eta_Y\p{\cap}(y\p{\cup}f^*(x)))=$
\hfill by (\ref{cap=cup})\\
$=f_!((\eta_Y\p{\cap}y)\p{\cap}f^*(x)) =$\hfill proj. form.
for $\p{\cap}$\\ $=f_!(\eta_Y\p{\cap}y)\p{\cap}x=$\hfill
by (\ref{cov})\\  $=(\eta_X\p{\cap}f_*(y))\p{\cap}x=$
\hfill by (\ref{cap=cup})\\
$=\eta_X\p{\cap}(f_*(y)\p{\cup}x).$} \E{center}

Because of the lack of symmetry of the projection formula
stated above we need to assume the following
``projection formula with supports''.
\B{description}
\item[$\forall$5] For $f:Y\to X$, $Z$ and $\rho$ as above, the
following:
$$\B{array}{c}
\hspace{20pt} H^{p}_{f^{-1}(Z)}(Y,\p{r})\tensor
H^{q}(Y,\p{s})\\ \p{id\tensor
f^{*}}\nearrow \ \hspace{40pt} \searrow
\p{\cup}  \\ H^{p}_{f^{-1}(Z)}(Y,\p{r})\tensor
H^{q}(X,\p{s})\hspace{60pt}
H^{p+q}_{f^{-1}(Z)}(Y,\p{r+s})\\  \p{f_{*}\tensor
id}\downarrow \hspace{121pt}\downarrow \p{f_{*}} \\
H^{p+2\rho}_{Z}(X,\p{r+\rho})\tensor
H^{q}(X,\p{s})\hspace{10pt} \longby{\cup}
\hspace{10pt} H^{p+q+2\rho}_{Z}(X,\p{r+s+\rho})
\E{array}$$
commutes, i.e. we have the equation
$f_*(y\p{\cup}f^*(x))=f_*(y)\p{\cup}x$ for $x\in
H^{q}(X,\p{s})$ and $y\in
H^{p}_{f^{-1}(Z)}(Y,\p{r})$
\E{description}

As a matter of fact, in order to obtain the projection
formula (\ref{chpr}) we just need the following apparently
weaker but almost equivalent form of $\forall$5.
\B{description}
\item[$\forall$5'] Let $i:Z\hookrightarrow X$ be a smooth pair of
pure codimension $c$. Then:
$$\B{array}{c}
\hspace{20pt}
H^{p}(Z,\p{r})\tensor H^{q}(Z,\p{s})\\
\p{id\tensor i^{*}}\nearrow \ \hspace{40pt} \searrow
\p{\cup}  \\ H^{p}(Z,\p{r})\tensor
H^{q}(X,\p{s})\hspace{60pt}
H^{p+q}(Z,\p{r+s})\\  \p{i_{*}\tensor
id}\downarrow\wr \hspace{121pt}\wr\downarrow
\p{i_{*}} \\ H^{p+2c}_{Z}(X,\p{r+c})\tensor
H^{q}(X,\p{s})\hspace{10pt} \longby{\cup}  \hspace{10pt}
H^{p+q+2c}_{Z}(X,\p{r+s+c})
\E{array}$$
commutes.
\E{description}

By the way $\forall$5 implies $\forall$5'. (Convention: the purity
isomorphism $i_*$ is induced by the identity on $Z$.)

\B{schol} Let $i: Z\hookrightarrow X$ be a smooth pair.
Then the following square

$$\begin{array}{ccc}
H^{p}_Z(X,\p{r})\tensor H^{q}(X,\p{s}) &\longby{\cup}
& H^{p+q}_Z(X,\p{r+s})\\
\p{\eta_X\cap -\tensor i^*}\downarrow & &\downarrow
\wr\ \p{\eta_X\cap -}\\
 H_{2d-p}(Z,\p{d-r})\tensor H^{q}(Z,\p{s})  &
\longby{\cap}&  H_{2d-p-q}(Z,\p{d-r-s})
\end{array}
$$
commutes, i.e. we have the following formula:
\B{equation}\label{rest}
\eta_X \p{\cap}(z\p{\cup}x) = (\eta_X
\p{\cap}z)\p{\cap}i^*(x)
\E{equation}
for $z\in H^{p}_Z(X,\p{r})$ and $x\in H^{q}(X,\p{s})$,
{\rm if and only if $\forall$5'} holds.
\E{schol}
\B{proof} Let $c$ be the codimension of $Z\subset X$ and
$d=$ dim $X$.\\ Let assume that $\forall$5' holds.  Since
we do have the purity isomorphism
$i_*:H^{p-2c}(Z,\p{r-c})\by{\cong} H^{p}_Z(X,\p{r})$ there
is an element $\zeta\in H^{p-2c}(Z,\p{r-c})$ such that
$i_*(\zeta)=z$. Because of $\forall$5' the equation (\ref{rest})
is obtained by showing the following equality
\B{equation}\label{zeta}
\eta_X\p{\cap}i_*(\zeta\p{\cup}i^*(x))=
(\eta_X\p{\cap}i_*(\zeta))\p{\cap}i^*(x) \E{equation}
Note: $\eta_Z\p{\cap}- :H^{p-2c}(Z,\p{r-c})
\by{\cong}  H_{2d-p}(Z,\p{d-r})$ and $\eta_Z\p{\cap}- =
\eta_X\p{\cap}i_*(-)$ by the definition of $i_*$; thus
the right-hand side in the equation (\ref{zeta}) above
becomes
\B{equation}\label{eta}
(\eta_Z\p{\cap}\zeta)\p{\cap}i^*(x)
\E{equation}
and the left-hand side becomes
$$
\eta_Z\p{\cap}(\zeta\p{\cup}i^*(x))
$$
which by $\forall$4 (see the equation (\ref{cap=cup})) is exactly
the same of (\ref{eta}).\\ Conversely, if (\ref{rest})
holds we then have
\B{center}
\parbox{3in}{$\eta_X\p{\cap}(i_*(\zeta)\p{\cup}x)=$\hfill
by (\ref{rest})\\
$=(\eta_X\p{\cap}i_*(\zeta))\p{\cap}i^*(x)=$\hfill
by def. of $i_*$\\
$=(\eta_Z\p{\cap}\zeta)\p{\cap}i^*(x)=$\hfill
by (\ref{cap=cup})\\
$=\eta_Z\p{\cap}(\zeta\p{\cup}i^*(x))=$\hfill
by def. of $i_*$\\
$=\eta_X\p{\cap}i_*(\zeta\p{\cup}i^*(x)).$}
\E{center}
Since $\eta_X\p{\cap}-$ is an isomorphism we can erase
it from the left of the resulting equation.
 \E{proof}
\B{rmk} As we will see below (cf. Lemma \ref{P5}): the
formula (\ref{rest}), by restriction and the projection
formula w.r.t. cap-product, allow us to deduce $\forall$5
for $X$, $Y$, $Z$ and $f^{-1}(Z)$ smooth. Moreover, by
assuming the formula (\ref{rest}) holds for $Z$ possibly
singular we can as well obtain $\forall$5.
\E{rmk}
\B{defi} Let $(H^*,H_*)$ be a Poincar\'e duality theory
with supports and let assume that $H^*$ has a cup-product
(so that $\forall$1--$\forall$3 are satisfied). We will
say that $(H^*,H_*)$ is {\it multiplicative}\, if the
axioms $\forall$4 and $\forall$5' (or the strong form
$\forall$5) are satisfied.
 \E{defi}

\subsection{$\cH$-cup products}\label{H-cup}

Let $H^*$ be a twisted cohomology theory with
supports on $\cat{V}_k$. Let $X\in
 \cat{V}_k$ be equidimensional but possibly singular.
Applying the exact couple method to the exact sequence
(\ref{limloc})
 $$
 H^i_{Z^{p+1}(X)}(X,\p{j})\to
H^i_{Z^p(X)}(X,\p{j})\to \coprod_{x\in X^p}^{}
H^{i}_x(X,\p{j}) \to H^{i+1}_{Z^{p+1}(X)}(X,\p{j})
$$
where: $H^{i}_x(X,\p{j}) \df \limdir{\pp{U\subset X}}
H^i_{\overline{\{x\}}\cap U}(U,\p{j})$, we do get the
coniveau spectral sequence  $$E^{p,q}_1 =\coprod_{x\in
X^p}^{} H^{q+p}_x(X,\p{\cdot}) \Rightarrow
H^{p+q}(X,\cdot)$$ Let denote ${\rm R}_q^{\bul}(X)\p{(r)}$
the corresponding Gersten type complexes $E^{\bul,q}_1$.

\B{prop} \label{pairing}Let $H^*$ be a twisted cohomology
theory with supports and cup-product on $\cat{V}_k$.
For $X\in\cat{V}_k$ there is a pairing of complexes
$$ {\rm R}_q^{\bul}(X)\p{(r)}\tensor H^n(X,\p{s})\to
{\rm R}_{q+n}^{\bul}(X)\p{(r+s)}$$
contravariant w.r.t. flat maps.
\E{prop}
\B{proof}
Let $Z\subseteq T\subseteq X$ with
$Z\in Z^{p+1}(X)$ and $T\in Z^p(X)$ and let $U=X-Z$; thus
by restriction to $U$ and cup-product we do have a pairing
associated to such pairs $Z\subseteq T$:
$$H_{T\cap U}^i(U,\p{r})\tensor H^n(X,\p{s}) \to
H_{T\cap U}^{i+n}(U,\p{r+s}) $$
i.e. $t\p{\tensor}x \leadsto t\p{\cup}j^*(x)$ where $j:
U\hookrightarrow X$. By taking the direct limit over
such pairs (this makes sense because of
$\forall$1--$\forall$2)  we do have a pairing
$$\coprod_{x\in X^p}^{} H^{q+p}_x(X,\p{r})\tensor
H^n(X,\p{s}) \to \coprod_{x\in X^p}^{}
H^{q+n+p}_x(X,\p{r+s})$$ In order to check compatibility
with the differentials of ${\rm R}_q^{\bul}(X)\p{(r)}$,
because of $\forall$2 we have a pairing
$H^{i}_{Z^{p}(X)}(X,\p{j})\tensor H^n(X,\p{s})\to
H^{i+n}_{Z^{p}(X)}(X,\p{j+s})$ and, by construction, the
differential is the composition of $${\rm
R}_q^{p}(X)\p{(r)}\to
H^{q+p+1}_{Z^{p+1}(X)}(X,\p{r})\to{\rm
R}_q^{p+1}(X)\p{(r)}$$ we can argue as in the proof of
the Proposition~\ref{hpairing} via $\forall$3 and limit
arguments.
\E{proof}

\B{defi} For $X\in\cat{V}_k$, by taking
associated sheaves for the Zariski topology of the pairing
above, we get a cap-product pairing
$${\cR}_q^{\bul}(X)\p{(r)}\tensor
\cH^p_X\p{(s)}\to {\cR}_{p+q}^{\bul}(X)\p{(r+s)}$$ By
sheafifying the cup-product we do have a product
$$\p{\cup}_{\cH}:\cH^p_X\p{(r)} {\tensor} \cH^q_X\p{(s)}\to
\cH^{p+q}_X\p{(r+s)}$$  which we call {\it
$\cH$-cup-product} on $X$.
\E{defi}

\B{rmk} We list several expected compabilities.\\
\B{enumerate}
\item The $\cH$-products above
are compatible via the canonical augmentations
$\cH^p_X\p{(r)}\to {\cR}_p^{\bul}(X)\p{(r)}.$  But, if $X$
is singular the augmentations are not quasi-isomorphisms.
\item Let suppose the existence of an external
product
 $$\p{\times}:
H_Z^p(X,\p{r})\tensor H_T^q(Y,\p{s})\to
H^{p+q}_{Z\times T}(X\times Y,\p{r+s})$$
functorial on $\cat{V}^2_k$ i.e. we have
the equation
\B{equation}(f\p{\times} g)^*(x\p{\times}y)
= f^*(x)\p{\times}g^*(y)
 \E{equation}
for $f$ and $g$ maps of pairs. Thus, by composing with
the diagonal $\Delta :(X,Z\p{\cap}T) \to (X\times X,
Z\times T)$, we obtain a cup-product satisfying the axiom
$\forall$1.\\ In this case the pairing defined in the
Proposition~\ref{pairing} can be obtained as follows
$$H_{T\cap U}^i(U,\p{r})\tensor H^n(X,\p{s})
\ni t\p{\tensor}x \leadsto \Delta^*(1\p{\times}j)^*(t
\p{\times}x)$$
where $j: U\hookrightarrow X$.
\item On a smooth variety  $X$, after $\forall$4, we have that the
pairing defined in Proposition~\ref{pairing} is Poincar\'e
dual of the pairing defined as follows:  let $j:
U\hookrightarrow X$, if $\xi\in H_*(X,\cdot)$ and $t\in
H^*_{\overline{\{p\}}\cap U}(U,\cdot)$ then
$\xi\p{\tensor}t \leadsto j^*(x)\p{\cap}t$. In fact:
$\xi=\eta_X\p{\cap}x$ for some  $x\in H^*(X,\cdot)$ and we
then have:
\B{center}
\parbox{3in}{
$\eta_U\p{\cap}(j^{*}(x)\p{\cup}t)=$\hfill by $\forall$4\\
$=(\eta_U\p{\cap}j^{*}(x))\p{\cap}t=$\hfill
\cite[1.3.2-4]{BO}\\
$=j^{*}(\eta_X\p{\cap}x)\p{\cap}t).$}
\E{center}
\item  On a smooth variety  $X$, after the equation
(\ref{rest}), we have that the pairing defined in the
Proposition~\ref{pairing} can be obtained by restriction
as follows. Let  $t\in H^*_{\overline{\{p\}}\cap
U}(U,\cdot)$ and  $x\in H^*(X,\cdot)$ where $X$ is smooth
and  $\overline{\{p\}}\cap U \subset U$ is a smooth pair.
Thus $\eta_U\p{\cap}(t\p{\cup}j^*(x))$ is equal to
$(\eta_U\p{\cap}t)\p{\cap}i^*(x)$
where $i: \overline{\{p\}}\cap U \hookrightarrow X$.
\E{enumerate}\E{rmk}

\subsection{$\cH$-Gysin maps}\label{H-Gys}

Let $(H^*,H_*)$ be a Poincar\'e duality theory with
supports. For $X\in\cat{V}_k$ we have a niveau spectral
sequence $E_{a,b}^1 \Rightarrow H_{a+b}(X,\cdot)$
which is covariant for proper morphisms. For $k$  perfect
and $X$ smooth equidimensional, $d=$ dim $X$,  by local
purity, we do have that $H^{q+p}_x(X,\p{r})\cong
H^{q-p}(x,\p{r-p})$ if $x\in X^p$ and isomorphisms

 $$E_{d-p,d-q}^1\p{(d-r)}
=\coprod_{x\in X_{d-p}}^{} H_{2d-p-q}(x,\p{d-r}) \cong
\coprod_{x\in X^p}^{}
H^{q-p}(x,\p{r})\cong E^{p,q}_1\p{(r)}$$

\B{lemma}\label{co=ni}If $X\in \cat{V}_k$ is smooth of pure
dimension $d$ then $E_{d-p,d-q}^1\p{(d-r)}\cong
E^{p,q}_1\p{(r)}$ is an isomorphism of spectral sequences
which is natural w.r.t. \'etale maps.
 \E{lemma}
\B{proof} This is a consequence of the above once we have
identified (via Poincar\'e duality) the long exact sequence
of cohomology with supports with the corresponding long
exact sequence of homology groups.
\E{proof}

Let $f:Y\to X$ be a proper morphism between
$k$-algebraic schemes where dim~$X$ = $\delta$ and dim
$Y$ = $d$. Let $\rho =\delta - d$. Since  $E^1\p{(\cdot)}$
is covariant w.r.t. proper maps we do have a map of niveau
spectral sequences  $$E_{d-p,d-q}^1\p{(d-r)}(Y) \to
E_{d-p,d-q}^1\p{(d-r)}(X)$$ If $X$ and $Y$ are {\it
smooth}\, and equidimensionals, then, via the Lemma
\ref{co=ni}, we get a map of coniveau spectral sequences
as follows: $$E^{p,q}_1\p{(r)}(Y)\cong
E_{d-p,d-q}^1\p{(d-r)}(Y) \to
E_{d-p,d-q}^1\p{(d-r)}(X)\cong E^{p+\rho ,q+\rho
}_1\p{(r+\rho )}(X)$$

\B{defi} For $f:Y\to X$ as above we will call the induced
map of complexes ${\rm R}_q^{\bul}(Y)\p{(r)}\to {\rm
R}_{q+\rho}^{\bul}(X)\p{(r+\rho)}[\rho]$ the {\it global
Gysin map}. By taking associated sheaves we have the
{\it local} Gysin map
$$f_{\flat}: f_*{\cR}_q^{\bul}(Y)\p{(r)}\to
{\cR}_{q+\rho}^{\bul}(X)\p{(r+\rho)}[\rho].$$ In the
derived category we do have the $\cH$-Gysin map
$${\bf R}f_{\flat}: {\bf R}f_*\cH^q_Y\p{(r)} \to
\cH^{q+\rho}_X\p{(r+\rho)}[\rho].$$
 \E{defi}

\B{rmk} For $i: Z\hookrightarrow X$ a (smooth) pair of
pure codimension $c$ we do have the isomorphism (cf.
\S\ref{cad1})
$$i_Z: i_*\cR_{q,Z}^{\bul}\p{(r)}
\by{\simeq}{\bf \Gamma}_Z\cR_{q+c,X}^{\bul}\p{(r+c)}[c].$$
Thus, as it is easily seen, the local Gysin map
$i_{\flat}$ is obtained by composition of $i_Z$ with the
canonical map ${\bf\Gamma}_Z\cR_{q+c,X}^{\bul}\p{(r+c)}[c]
\hookrightarrow \cR_{q+c,X}^{\bul}\p{(r+c)}[c]$ (see the
proof of Scholium~\ref{Gysin}).
 \E{rmk}

\subsection{Projection formula}

Let $f:Y\to X$ be a proper morphism between smooth
equidimensional algebraic schemes over a perfect field
$k$. Let dim $X$ = $\delta$ and dim $Y$ = $d$. Let $\rho
=\delta - d$.

\B{prop} For $f:Y\to X$ as above we have the following
commutative diagram:

\B{equation}
\B{array}{c}\label{chpr}
\hspace{10pt} f_*\cR^{\bul}_{q,X}\p{(r)}\tensor
f_*\cH^{p}\p{(s)}\\ \p{id\tensor f^{\sharp}}\nearrow \
\hspace{40pt} \searrow \p{f_*\cap_{{\cH}}}  \\
f_*\cR^{\bul}_{q,X}\p{(r)}\tensor
\cH^{p}\p{(s)}\hspace{60pt}
f_*\cR^{\bul}_{p+q,X}\p{(r+s)}\\  \p{f_{\flat}\tensor
id}\downarrow \hspace{105pt}\downarrow \p{f_{\flat}} \\
\cR^{\bul}_{q+\rho ,X}\p{(r+\rho )}[\rho ]\tensor
\cH^{p}\p{(s)} \hspace{10pt} \longby{\cap_{\cH}}
\hspace{10pt}\cR^{\bul}_{p+q+\rho,X}\p{(r+s+\rho)}[\rho]
\E{array}
\E{equation} \E{prop}

\B{proof} The commutative diagram above will be otbained
from the following:
$$
\B{array}{c}
\hspace{20pt} {\rm R}_q^{\bul}(f^{-1}(U))\p{(r)}\tensor
H^{p}(f^{-1}(U),\p{s})\\ \p{id\tensor f^{*}}\nearrow \
\hspace{40pt} \searrow \p{\cap}  \\
{\rm R}_q^{\bul}(f^{-1}(U))\p{(r)}\tensor
H^{p}(U,\p{s})\hspace{60pt}
{\rm R}_{p+q}^{\bul}(f^{-1}(U))\p{(r+s)}\\  \p{f_{\flat}
\tensor id}\downarrow \hspace{105pt}\downarrow
\p{f_{\flat}} \\ {\rm
R}_{q+\rho}^{\bul}(U)\p{(r+\rho)}[\rho]\tensor
H^{p}(U,\p{s}) \hspace{10pt} \longby{\cap}  \hspace{10pt}
{\rm R}_{p+q+\rho}^{\bul}(U)\p{(r+s+\rho)}[\rho]
\E{array} $$
where $U\subset X$ is any Zariski open subset of $X$, by
taking associated sheaves on $X_{Zar}$.\\
Moreover it sufficies to prove the case of $U=X$.\\
Let $\overline{\{y\}}\subset Y$ such that
$y\in Y^{c}$ and $f(y)\in X^{c+\rho}$. The Gysin map (cf.
\S\ref{H-Gys}) $$f_{\flat,y}: \limdir{\pp{V\subset
Y}}H^*_{\overline{\{y\}}\cap V}(V) \to \limdir{\pp{U\subset
X}}H^*_{\overline{\{f(y)\}}\cap U}(U)$$
is the Poincar\'e dual of  $f_{!,y}:
\limdir{\pp{V\subset Y}}H_*({\overline{\{y\}}\cap V}) \to
\limdir{\pp{U\subset X}}H_*({\overline{\{f(y)\}}\cap U})$
(see the proof of (\ref{hpr})).
By the definition of the pairing in Proposition
\ref{pairing} we are left to show that the following
diagram

\B{equation}\B{array}{c} \label{open}
\ \hspace{20pt}
H_{\overline{\{y\}}\cap V}^p(V,\cdot)\tensor H^q(V,\cdot)
\\ \p{id\tensor f^*}\nearrow \hspace{50pt}
\searrow\p{\cup}\\  H_{\overline{\{y\}}\cap
V}^{p}(V,\cdot)\tensor H^q(U,\cdot)\hspace{40pt}
H_{\overline{\{y\}}\cap V}^{p+q}(V,\cdot) \\
\p{f_{\flat}\tensor id}\downarrow \hspace{105pt} \downarrow
\p{f_{\flat}} \\ H_{\overline{\{f(y)\}}\cap U}^{p+2\rho}(U,\cdot
+\rho)\tensor H^q(U,\cdot) \hspace{10pt} \longby{\cup}
\hspace{10pt} H_{\overline{\{f(y)\}}\cap
U}^{p+q+2\rho}(U,\cdot +\rho)  \E{array}
\E{equation}
commutes where $f^{-1}(U) = V$. By shrinking the  open
sets involved we may assume that $\overline{\{y\}}\cap V
\subset V$ and  $\overline{\{f(y)\}}\cap U \subset U$ are
smooth pairs. Since $f(\overline{\{y\}}\cap f^{-1}(U))=
\overline{\{f(y)\}}\cap U$ the diagram (\ref{open})
commutes because of the following Lemma.
\E{proof}

\B{lemma}\label{P5} Let $f:Y\to X$ and $\rho$ as above.
Let $T$ be a closed subset of $Y$, $f(T)=Z$ and let assume
that  $T\hookrightarrow Y$ and $Z \hookrightarrow X$ are
smooth pairs. Then the following diagram

\B{equation}\B{array}{c}
\ \hspace{20pt}
H_{T}^p(Y,\p{r})\tensor H^q(Y,\p{s})
\\ \p{id\tensor f^*}\nearrow \hspace{50pt}
\searrow\p{\cup}\\  H_{T}^{p}(Y,\p{r})\tensor
H^q(X,\p{s})\hspace{40pt} H_{T}^{p+q}(Y,\p{r+s}) \\
\p{f_{\flat}\tensor id}\downarrow \hspace{105pt} \downarrow
\p{f_{\flat}} \\ H_{Z}^{p+2\rho}(X,\p{r
+\rho})\tensor H^q(X,\p{s}) \hspace{10pt} \longby{\cup}
\hspace{10pt} H_{Z}^{p+q+2\rho}(X,\p{r+s+\rho})  \E{array}
\E{equation}
commutes.
\E{lemma}
\B{proof} We will give two proofs.\\
{\it First proof.}\,  Let assume the axiom $\forall$5 (in
which case we do not need the smoothness of $T$ and $Z$).
Let  $i: T \hookrightarrow f^{-1}(Z)$ so that $f\mid_{T} =f
\p{\circ} i$ hence $f\mid_{T,!} =f_!\p{\circ} i_!$ and
$f_{\flat} =f_*\p{\circ} i_{\diamond}$ where
$i_{\diamond}: H^*_T(Y,\cdot) \to H^*_{f^{-1}(Z)}(Y\cdot)$
is the canonical map. Thus, for $y\in H^*_T(Y,\cdot)$ and
$x\in H^*(X,\cdot)$
 \B{center}
\parbox{3in}{$f_{\flat}(y\p{\cup}f^*(x))=$ \\
$=f_{*}\p{\circ}i_{\diamond}(y\p{\cup}f^*(x))=$\hfill by
$\forall$2\\  $=f_{*}(i_{\diamond}(y)\p{\cup}f^*(x))=$\hfill by
$\forall$5\\ $=f_{*}(i_{\diamond}(y))\p{\cup}x)=$\\
$=f_{\flat}(y)\p{\cup}x.$}
\E{center}

\noindent{\it Second proof.}\,  Just assuming the axiom
$\forall$5' we can prove the Lemma as follows.
Let $k:T\hookrightarrow Y$,  $i:Z\hookrightarrow
X$ and $f\mid_T: T \to Z$. Thus we have:
$k^*f^*=(f\mid_{T})^*i^*$.  Let $y\in H_{T}^{p}(Y,\p{r})$
and $x\in H^q(X,\p{s})$. We have:

\B{center}
\parbox{5in}{$\eta_X\p{\cap}f_{\flat}(y\p{\cup}f^*(x))=$
\hfill by def. of $f_{\flat}$\\
$=(f\mid_{T})_!(\eta_Y\p{\cap}(y\p{\cup}f^*(x)))=$\hfill
by (\ref{rest})\\
$=(f\mid_{T})_!((\eta_Y\p{\cap}y)\p{\cap}k^*f^*(x))=$\ \\
$=(f\mid_{T})_!((\eta_Y\p{\cap}y)\p{\cap}
(f\mid_{T})^*i^*(x))=$\hfill by \cite[1.3.3]{BO}\\
$=(f\mid_{T})_!(\eta_Y\p{\cap}y)\p{\cap}i^*(x)=$\hfill by
def. of $f_{\flat}$\\
$=(\eta_X\p{\cap}f_{\flat}(y))\p{\cap}i^*(x)=$\hfill  by
(\ref{rest})\\ $=\eta_X\p{\cap}(f_{\flat}(y)\p{\cup}x)$.}
\E{center}
Since $\eta_X\p{\cap}-$ is an isomorphism we conclude.
\E{proof}

\B{cor} For $i: Z\into X$ a smooth pair of pure
codimension $c$ the following diagram

\B{equation}
\B{array}{c}\label{ipr}
\hspace{10pt} i_*\cR^{\bul}_{q,Z}\p{(r)}\tensor
i_*\cH^{p}\p{(s)}\\ \p{id\tensor i^{\sharp}}\nearrow \
\hspace{40pt} \searrow \p{i_*\cap_{{\cH}}}  \\
i_*\cR^{\bul}_{q,Z}\p{(r)}\tensor
\cH^{p}\p{(s)}\hspace{60pt}
i_*\cR^{\bul}_{p+q,Z}\p{(r+s)}\\  \p{i_{Z}\tensor
id}\downarrow\wr \hspace{105pt}\wr\downarrow \p{i_{Z}} \\
{\bf \Gamma}_Z\cR^{\bul}_{q+c,X}\p{(r+c)}[c]\tensor
\cH^{p}\p{(s)} \hspace{10pt} \longby{\cap_{\cH}}
\hspace{10pt}{\bf \Gamma}_Z\cR^{\bul}_{p+q+c,X}\p{(r+s+c)}[c]
\E{array}
\E{equation}
commutes.
 \E{cor}
\B{proof} This follows by the factorisation of $i_{\flat}$
(cf. the Remark at the end of \S\ref{H-Gys}) and the
Proposition above.
 \E{proof}

\subsection{$\cH$-cohomology ring}\label{H-ring}

Let $(H^*,H_*)$ be a multiplicative Poincar\'e duality
theory with supports. Let suppose that our cohomology
theory $H^*$ takes values in a fixed category of
$\La$-modules where $\La = H^0(k,\p{0})$ is a commutative
ring with 1; we assume that the bigraded $\La$-module
$\bigoplus_{q,r}^{} H^q(X,\p{r})$ has a $\La$-algebra
structure via the cup-product pairing e.g. a canonical
isomorphism of rings $H^0(X,\p{0})\cong \La$ if $X$ is
irreducible.\\ Let $$A(X) \df \bigoplus_{p,q,r}^{}
H^p(X,\cH^q\p{(r)}).$$ Then $X \leadsto A(X)$ is a
contravariant functor on $\cat{V}_k$. If $f:Y \to X$ is a
proper map of relative dimension $\rho$ then the
$\cH$-Gysin maps ${\bf R}f_*\cH^q\p{(r)}\to
\cH^{q+\rho}\p{(r+\rho)}[\rho]$ induce direct image $A(Y)
\to A(X)$ (a map of degree $\rho$) so that $A$ is a
covariant functor w.r.t. proper maps of smooth varieties.
{}From the $\cH$-cup-product pairing by taking cohomology we
have an external pairing
$$\p{\times}:A(X)\otimes_{\La}A(Y) \to A(X\times Y)$$
which is associative and anticommutative (can be made
commutative by using the trick in \cite{GIN}). In
particular, let consider the functor   $$X \leadsto
\bigoplus_{p}^{}  H^p(X,\cH^p\p{(p)})\df A_{diag}(X).$$
If $\cH^0(\p{(0)})$ is identified with the flasque sheaf
$\coprod_{X^0}^{} \La$ (e.g. by assuming the `dimension
axiom')  we then have an augmentation $\varepsilon
:A_{diag}(X)\to\La$ where $X$ is irreducible and
$\varepsilon^0 :  H^0(X,\cH^0\p{(0)})\cong \La$  zero
otherwise.\\ Let denote $f_*$ and $f^*$ ``direct and
inverse'' images. We have the following formulas (cf.
\cite[I.1--I.9]{GI}): \B{eqnarray}
(f\p{\times}g)^*(-\p{\times}\cdot)&
=& f^*(-)\p{\times}g^*(\cdot)\nonumber \\
(f\p{\times}g)_*(-\p{\times}\cdot)
&=&f_*(-)\p{\times}g_*(\cdot)\nonumber \\
\varepsilon (f^*(\cdot))&=&\varepsilon (\cdot)
\nonumber \\  \varepsilon
(-\p{\times}\cdot)&=&\varepsilon (-)\varepsilon (\cdot)
 \E{eqnarray}
Furthermore, for $X={\rm Spec} k$ we have that
$\varepsilon : A_{diag}(k)\cong \La$: let $e$ be the
unique element such that $\varepsilon (e) =1$. We have
the formula: $e\p{\times}-=-\p{\times}e=-$.\\ Let
$\lambda :X\to {\rm Spec} k$ be the structural map and
let denote $\lambda^*(e)\df 1_X$. This is equal to
$\varepsilon^{-1}(1)$ on $X$ irreducible. For $x\in
A_{diag}(X)$ and $y\in A_{diag}(Y)$ we have
\B{eqnarray} \label{cross}
x\p{\times}1_Y&=&p_1^*(x)\nonumber\\
1_X\p{\times}y&=&p_2^*(y)
\E{eqnarray}
where $p_1$ and $p_2$ are the first and the second
projections of $X\times Y$ on its factors.\\
By composing the external product
$\p{\times}:A_{diag}(X)\otimes_{\La}A_{diag}(X) \to
A_{diag}(X\times X)$  with the diagonal $\Delta^*_X:
A_{diag}(X\times X)\to A_{diag}(X)$ we do get a product
$x\p{\tensor}x\prime \leadsto xx\prime$ in $A_{diag}(X)$
making it an associative anticommutative algebra with
identity $1_X$. The homomorphism $\varepsilon :
A_{diag}(X)\to \La$ is an homomorphism of unitary
$\La$-algebras. For $f:X\to Y$ the map $f^*:
A_{diag}(Y)\to A_{diag}(X)$ is a homomorphism of
$\La$-algebras. The external product is a homomorphism
of augmented $\La$-algebras. This last fact, via the
equations (\ref{cross}), give us the formula
\B{equation}
x\p{\times}y=p_1^*(x)p_2^*(y)
\E{equation}
for $x\in A_{diag}(X)$ and $y\in A_{diag}(Y)$.\\
For  $f:X\to Y$ a proper map of smooth varieties
over $k$, $A_{diag}$ satisfies the `projection formula' as
a consequence of the projection formula (\ref{chpr}).
Furthermore, if $f$ is surjective of relative dimension
$\rho$ over $Y$ irreducible we have a canonical map
$$\int_{X/Y}: H^{-\rho}(X,\cH^{-\rho}\p{(-\rho)})\to \La$$
and its extension by zero $\int : A_{diag}(X)\to\La$,
both defined by composition of the $\cH$-Gysin map $f_*$
and the augmentation $\varepsilon$. In particular
$$\int_{X/k}: H^{d}(X,\cH^{d}\p{(d)})\to \La$$
for any $X$ proper smooth $d$-dimensional variety; for
any map $f$ between $X$ and $Y$ proper smooth varieties we
have $$\int_{Y/k}f_* = \int_{X/k}$$ Finally, if the proper
map has a section $fs=1$ then $\int$ is a surjection; this
is the case of $X/k$ having a $k$-rational point.

\section{Intersection theory}

Since we are going to deal with Poincar\'e duality
theories which are `sophisticated' {\it and}\,
`multiplicatives' we need to arrange the axioms in order
to be not redundant. This arrangement will yields the
notion of `duality theory {\it appropriate}\, for algebraic
cycles' or for short `appropriate duality theory'. We will
show that the $\cH$-cohomology rings associated with such
a theory reproduce the classical intersection rings.
Roughly speaking, the denomination `appropriate duality'
is the corresponding cohomological version of `relation
d'\'equivalence ad\'equate' introduced by P.Samuel (see
\cite{SAM}) for algebraic cycles.

\subsection{Axiomatic menuet}

Let $H^*$ be a cohomology theory and let $H_*$ be a
homology theory (as defined by \cite[1.1 and 1.2]{BO}).
Let assume that the pair $(H^*,H_*)$ yields a sophisticated
Poincar\'e duality which satisfies the dimension axiom;
furthermore we assume the existence of an associative
anticommutative functorial cup product pairing $$
H_Y^p(X,\p{r})\tensor H_Z^q(X,\p{s})\to H^{p+q}_{Y\cap
Z}(X,\p{r+s})$$ where $\La = H^0(k,\p{0})$ is a commutative
ring with 1 and the bigraded $\La$-module
$\bigoplus_{q,r}^{} H^q(X,\p{r})$ has a $\La$-algebra
structure via the cup-product pairing.\\[3pt]

\B{defi} An {\it appropriate}\, duality theory is a pair
$(H^*,H_*)$ as above such that the sophisticated cap
product is compatible with the cup product via Poincar\`e
duality i.e. the following diagram, where $q+j=2d, s+n=d$
and $X$ is smooth
\B{displaymath}\B{array}{ccc}
 H^q_Y(X,\p{s})\tensor H_Z^p(X,\p{r})  & \longby{\cup} &
H^{p+q}_{Y\cap Z}(X,\p{r+s})\\  \p{\eta_X\cap -\tensor id}
\downarrow & &  \downarrow\p{\eta_X\cap -} \\H_j(Y,\p{n})
\tensor H_Z^p(X,\p{r}) & \longby{\cap} &
H_{j-p}(Y\cap Z,\p{n-r})
\E{array}\E{displaymath}
commutes. In particular, the fundamental class $\eta_X\in
H_{2d}(X\p{d})$ corresponds to the unit $1\in\La\cong
H^0(X,\p{0})$ in the $\La$-algebra structure.
\E{defi} \hfill\\[4pt]
Let $(H^*,H_*)$ be an appropriate duality theory on
$\cat{V}_k$. Then, by adopting the same notation of
\S\ref{algcyc},
 $$\cat{V}_k\ni Y\leadsto
C_{n,m}(Y;\La\p{(s)})\df\H^{-n}(Y,{\cQ}^m_{\bul}\p{(s)})$$
is a covariant functor for proper morphisms in
$\cat{V}_k$. It is a presheaf for the \'etale
topology (or just for the Zariski topology, depending with
the homology theory). On the other hand we have a
contravariant functor  $$(X,Y)\leadsto
H^p_Y(X,\cH^q\p{(r)})$$ which yields the $\cH$-cohomology
ring with the properties stated in \S\ref{H-ring}. Indeed,
on a smooth scheme $X$ of pure dimension $d$, these two
functors are related via the duality isomorphism
$${\cQ}^{d-q}_{-\bul}\p{(d-r)}[d]\cong
{\cR}_q^{\bul}\p{(r)}$$ and this isomorphism is compatible
with the $\cH$-cap and $\cH$-cup products (cf.
\S\ref{H-cap} and \S\ref{H-cup}); by construction this
duality isomorphism identifies Gysin maps (cf.
\S\ref{H-Gys}) and projection formulas (cf. (\ref{hpr})
with (\ref{chpr})). Thus we do have a canonical
cap-product associated with pairs $(X,Y)$ and $(X,Z)$
  $$C_{n,m}(Y;\La\p{(s)})\tensor H^p_Z(X,\cH^q\p{(r)}) \to
C_{n-p,m-q}(Y\cap Z;\La\p{(s-r)})$$
and a corresponding projection formula.
We have indeed a canonical ``trace map'' on $X$
irreducible
 $${\cQ}^{d}_{-\bul}\p{(d)}[d]\by{\sim} \La $$
yielding a global section $[X]\in C_d(X;\La)$; by capping
with this ``fundamental class'' $[X]$ we get the
quasi-isomorphism $$\cH^q\p{(r)}\by{\cap [X]}
{\cQ}^{d-q}_{-\bul}\p{(d-r)}[d]$$
By taking hypercohomology with support on $Z$ (a closed
equidimensional subscheme) we do get the ``duality''
isomorphism
\B{equation}
\p{\cap}[X]: H^p_Z(X,\cH^q\p{(r)}) \by{\sim}
C_{d-p,d-q}(Z;\La\p{(d-r)}).
\E{equation}
Conversely: for $i:Z\into X$ we have a quasi-isomorphism
$$i_*{\cQ}^{d-q}_{-\bul}\p{(d-r)}[d]\by{\sim} {\bf
\Gamma}_Z\cH^q\p{(r)}$$ hence the canonical isomorphism
\B{equation}
\eta : C_{d-p,d-q}(Z;\La\p{(d-r)})\by{\sim}
H^p_Z(X,\cH^q\p{(r)})
\E{equation}
In particular:
\B{schol}\label{prs} We have a commutative diagram
\B{displaymath}
\B{array}{c}\hspace{10pt} C_{d-p}(Z;\La)\tensor
H^q_{Y\cap Z}(Z,\cH^q\p{(q)})\\ \p{id\tensor
i^*}\nearrow \ \hspace{40pt} \searrow
\p{H^*(\cap_{{\cH}})}  \\ C_{d-p}(Z;\La)\tensor
H^q_{Y}(X,\cH^q\p{(q)})\hspace{60pt}
C_{d-p-q}(Y\cap Z;\La)\\  \p{\eta\tensor
id}\downarrow \hspace{105pt}\downarrow \p{\eta}
\\ H^p_{Z}(X,\cH^p\p{(p)})\tensor
H^q_{Y}(X,\cH^q\p{(q)})\hspace{10pt}
\longby{H^*(\cup_{\cH})}  \hspace{10pt}H^{p+q}_{Y\cap
Z}(X,\cH^{p+q}\p{(p+q)})
  \E{array}\E{displaymath}
\E{schol}
\B{proof} This is a consequence of the commutative
diagram (\ref{hipr}) and the compatibilities between
$\cH$-products. \E{proof}

\subsection{$\cH$-cycle classes} \label{H-cycle}

We mantain the notations and the assumptions of the
previous section. Let $Z\subset X$ be a prime cycle of
dimension $d-c$. Then $$C_{d-c}(Z;\La) \df
C_{d-c,d-c}(Z;\La\p{(d-c)})=H_{2d-2c}(K(Z),\p{d-c})
\cong\La$$ by the `dimension axiom' and we do have a cycle
class
\B{equation}
\eta (Z) \in H^{c}_Z(X,\cH^{c}\p{(c)})
\E{equation}
where: $[Z]\in H_{2d-2c}(K(Z),\p{d-c})$ is obtained
by restriction of the fundamental class $\eta_Z\in
H_{2d-2c}(Z,\p{d-c})$ to the generic point and we have
$$\La\ni 1 \leadsto [Z] \leadsto \eta(1)\df \eta (Z)
\in   H^{c}_Z(X,\cH^{c}\p{(c)}).$$
Furthermore, by capping with the fundamental class $[X]$,
we find the formula  \B{equation}\label{inv} [Z]=\eta
(Z)\p{\cap}[X] \E{equation}
In particular the cycle class $\eta (Z)$ is independent
from the imbedding of $Z$ as a subvariety and it is
functorial w.r.t. \'etale maps.

\B{lemma}\label{etaext}
For $Y$ and $Z$ prime cycles of codimension $p$ and $q$
in $X$ smooth  we have
$$\eta(Y\p{\times}Z) =\eta (Y)
\p{\times}\eta (Z) \in H^{p+q}_{Y\times Z}(X\times
X,\cH^{p+q}\p{(p+q)}).$$ \E{lemma} \B{proof}
By standard sheaf theory the external product is obtained
by using flasque resolutions (see \cite[6.2.1]{GO}). Thus
via the canonical quasi-isomorphisms $\cH\p{()}
\cong\cQ_{\bul}\p{()}$ we do have a commutative diagram
\B{displaymath}\B{array}{ccc}
H_Y^p(X,\cH^p\p{(p)})\tensor H_Z^q(X,\cH^q\p{(q)}) &
\longby{\times} & H^{p+q}_{Y\times Z}(X\times
X,\cH^{p+q}\p{(p+q)})\\
 \downarrow\wr& &  \downarrow\wr \\
C_{d-p}(Y;\La)\tensor C_{d-q}(Z;\La) & \longby{\times} &
C_{d-p-q}(Y\times Z;\La)
\E{array}\E{displaymath}
To conclude one would see that the bottom arrow is in fact
the external product of cycles: this last claim is clear
because the external products are homomorphisms of
$\La$-algebras. (Note: $C_{dim?}(?;\La)\cong \La$ by the
dimension axiom).
\E{proof}

Let $\Delta :X\to X\times X$ be the diagonal embedding
and let $\Delta (X)$ be the diagonal cycle on $X\times X$.
Let $Y$ and $Z$ be prime cycles of codimension $p$ and $q$
in $X$ smooth such that $Y\cap Z$ is of pure codimension
$p+q$. For $d=$dim$X$ we then have
\B{displaymath}\B{array}{ccc}
H_{Y\cap Z}^{p+q}(X,\cH^{p+q}\p{(p+q)}) &
\longby{\Delta_{\bul}} & H^{p+q+d}_{(Y\times
Z)\cap \Delta (X)}(X\times X,\cH^{p+q+d}\p{(p+q+d)})\\
\p{\eta}\ \uparrow\wr& &  \uparrow\wr\ \p{\eta} \\
C_{d-p-q}(Y\cap Z;\La) & \longby{\Delta_*} &
C_{d-p-q}((Y\times Z)\cap \Delta (X);\La)
\E{array}\E{displaymath}
where $\Delta_{\bul}$ is obtained by making the diagram
commutative and $\Delta_*$ is induced by the isomorphism
$\Delta : Y\cap Z \to (Y\times Z)\cap \Delta (X)$. Thus
the formula
\B{equation}\label{etapoint}
\Delta_{\bul}(\eta (\ ))=\eta(\Delta_*(\ ))
\E{equation}

\B{lemma}
We have the following formula:
\B{equation}\label{etadiag}
\Delta_{\bul}\Delta^*(\eta(Y\times Z))=\eta(Y\times
Z)\eta (\Delta (X))
 \E{equation}
\E{lemma}
\B{proof}
The formula is a consequence of the projection formula
in the Scholium~\ref{prs} (cf. (\ref{hipr}) and
(\ref{ipr}) ) applied to the diagonal embedding $\Delta$ by
taking $\cH$-cohomology with supports on $Y\times Z$ and
compatibility of the $\cH$-products with the canonical
augmentations. \E{proof}

\subsection{Intersection of cycles}\label{int-cyc}

Let now assume that our cycle group $C_{*}(X;\La )$ has an
intersection product satisfying the classical
properties (cf. \cite{FU}): local nature of the
intersection multiplicity, normalization and reduction to
the diagonal for $X$ a smooth projective variety over a
field $k$.\\ Moreover, for a pair $(X,D)$ where
$X\in\cat{V}_k$ and $D$ is a Cartier divisor on $X$ we let
assume the existence of a homomorphism (cf.
Definition~\ref{line} in \S\ref{H-chern} below)
$$c\ell:H^1_D(X,\cO^*_X)\tensor \La\to
H^1_D(X,\cH^1\p{(1)})$$ such that  \B{description}
\item[{\it (i)}] $c\ell$ is a natural trasformation of
contravariant functors w.r.t. morphisms $f:X'\to X$ such
that $f^{-1}(D)$ is a divisor on $X'$;
\item[{\it (ii)}] $c\ell$ is compatible with the cap-products
in the sense that the following
\B{displaymath}\B{array}{ccc}
H_{D}^{1}(X,\cO^*_X)\tensor \La & \by{c\ell} &
H^1_D(X,\cH^1\p{(1)}) \\  \p{\cap [X]}\ \downarrow& &
\downarrow\ \p{\cap [X]} \\ CH_{d-1}(D;\La) & = &
C_{d-1}(D;\La)  \E{array}\E{displaymath} commutes, where
$d=$dim$X$.
\E{description}

\B{rmk} The map $\p{\cap [X]}: H_{D}^{1}(X,\cO^*_X)\to
CH_{d-1}(D)$ is given by the cap-product in
$\cK$-cohomology with the canonical cycle $[X]\in
CH_{d}(X)$; when applyed to the cycle class of the
Cartier divisor yields just the associated Weil divisor
(cf. \cite[\S2]{GIN}). By the way, if $X$ is non-singular
then $c\ell$ is an isomorphism.
 \E{rmk}

Thus we can prove the following key lemma.
\B{lemma}\label{intdiv} Let $X$ be smooth. Let $D$ be a
principle effective Cartier divisor and let $i: Z\into X$
be a closed integral subscheme of codimension $c$ in $X$
such that $Z\cap D$ is a divisor on $Z$. Then the
following  \B{displaymath}\B{array}{ccc}
 H^1_D(X,\cH^1\p{(1)}) &
\by{i^*} & H^1_{Z\cap D}(Z,\cH^1\p{(1)}) \\
\p{\eta}\ \uparrow& &  \downarrow\ \p{\cap [Z]} \\
C_{d-1}(D;\La) & \by{i^*} & C_{d-c-1}(Z\cap D;\La)
\E{array}\E{displaymath}
commutes i.e. we have the following formula
\B{equation}
D\bul Z = i^*\eta(D)\p{\cap} [Z]
\E{equation}
\E{lemma}
\B{proof} The claimed commutative diagram is obtained by
the corresponding one for the Picard groups (cf.
\cite[\S2]{GIN}). Let denote $\bar D\in
H_{D}^{1}(X,\cO^*_X)\tensor \La$ the canonical class of
the Cartier divisor: thus $\bar D\p{\cap} [X] = [D]\in
CH_{d-1}(D;\La)$; since $X$ is non-singular, $[D] \leadsto
\bar D$ under the isomorphism $CH_{d-1}(D;\La)\cong
H_{D}^{1}(X,\cO^*_X)\tensor\La$ and  $[D] \leadsto
\eta(D)$ under the isomorphism $CH_{d-1}(D;\La)\cong
H^1_D(X,\cH^1\p{(1)}).$\\ We then have
\B{center} \parbox{3in}{$i^*\eta(D)\p{\cap}
[Z] =$\hfill by {\it (ii)}\\ $=i^*c\ell(\bar D)\p{\cap} [Z]
=$\hfill by {\it (i)}\\ $=c\ell(i^*\bar D)\p{\cap} [Z]
=$\hfill by \cite[\S2]{GIN}\\ $=i^*(D) =$\hfill\\ $=D\bul
Z$} \E{center} where the last equality is just the
normalization property of the intersection theory.
\E{proof}
\B{teor}
With the above assumptions and notations, let $Y$ and $Z$
be prime cycles of codimension $p$ and $q$ on $X$ smooth
which intersect properly. Then $$ \eta (Y)\eta (Z) = \eta
(Y\bul Z) \in H^{p+q}_{Z\cap Y}(X,\cH^{p+q}\p{(p+q)})$$
\E{teor}
\B{proof} The proof is similar to that of the ``uniqueness
of the intersection theory'' and it consists of 3 steps.\\
{\it Step 1.\, (Intersection with divisors).} Let
$Y=D\into X$ be a principle Cartier divisor. Let $i:Z\into
X$. Then \B{center} \parbox{3in}{$\eta (D)\eta (Z)=$\hfill
by Scholium~\ref{prs}\\$\eta (i^*\eta (D) \p{\cap}
[Z])=$\hfill by Lemma~\ref{intdiv}\\$\eta (D\bul Z).$}
\E{center}
{\it Step 2.\, (Intersection with smooth subvarieties).}
Let assume $Y$ to be smooth. Since we can reduce to open
Zariski neighborhoods of the generic points of $Y\cap Z$
we may assume that $X$ is affine and $Y=V(f_1,\ldots
,f_p)$ where $\{f_1,\ldots ,f_p\}$ is a regular sequence.
Thus: $Y={\displaystyle\cap_{i=1}^{p} D_i}$ where $D_i
=V(f_i)$ and
\B{center} \parbox{3in}{$\eta (Y)\eta (Z) =$ \\
$=\eta (D_1\cdots D_p)\eta (Z)=$\\
$=\eta (D_1)\eta (D_2 \cdots D_p)\eta (Z)=\ldots$\\
$\ldots=\eta (D_1)\cdots \eta (D_{p-1})\eta (D_p\bul Z)=
\ldots$\\ $\ldots =\eta (Y\bul Z)$}
\E{center}
by iterative application of Step 1.\\
{\it Step 3.\, (Reduction to the intersection with the
diagonal).} We prove the general case as follows:
\B{center} \parbox{3in}{$\Delta_{\bul}(\eta (Y)\eta (Z))
=$\hfill by definition\\$=\Delta_{\bul}\Delta^*(\eta
(Y)\times \eta (Z))=$\hfill by Lemma~\ref{etaext}\\
$=\Delta_{\bul}\Delta^*(\eta (Y\times Z))=$\hfill by
the formula (\ref{etadiag})\\$=\eta (Y\times Z)\eta (\Delta
(X))=$\hfill by Step 2\\$=\eta (Y\times Z\bul \Delta
(X))=$\hfill int. with the diag.\\$=\eta
(\Delta_*(Y\bul Z)=$\hfill by the formula
(\ref{etapoint})\\ $=\Delta_{\bul}(\eta (Y\bul Z))$}
\E{center}
Since $\Delta_{\bul}$ is an isomorphism we conclude.
 \E{proof}

\B{cor} If $X$ is smooth of pure dimension $d$ then the
graded isomorphism $$\eta : \bigoplus C_{d-p}(X;\La) \cong
\bigoplus H^{p}(X,\cH^{p}\p{(p)})$$
is a $\La$-algebra isomorphism.
\E{cor}

\section{Chern classes and blow-ups}

Let $X$ be a variety i.e. $X\in \cat{V}_k$ reduced and
equidimensional over a  perfect field $k$, which admits a
closed imbedding in a smooth variety; such varieties
are usually called {\it imbeddable}. The existence of
$\cH$-cap-products grant us to construct Gysin maps for
the functor $C_{*}(-;\La)$ associated with such
imbeddings. By using the results from \S3--\S6 we construct
Chern classes in $\cH$-cohomologies. Furthermore, we are
able to obtain the nice decomposition formula for the
$\cH$-cohomology of blow-ups generalising the classical
one for Chow groups.

\subsection{Gysin maps for algebraic
cycles}

Let $(H^*,H_*)$ be an appropriate duality. We consider an
imbeddable variety $X$ with a fixed ambient smooth
variety $Y$. Let $i:X\into Y$ be a closed imbedding of
pure codimension $c$. Thus we have a $\cH$-cycle
class $\eta (X)\in H^c_X(Y,\cH^c\p{(c)})$ and the
corresponding Gysin maps
$$i^!:C_{n}(Y;\La)\to C_{n-c}(X;\La)$$
are defined as follows:
\B{equation}y\leadsto y\p{\cap}\eta (X)\df i^!(y)
\E{equation}
Thus $i^!(Y) = [X]$ because of $[Y]\p{\cap}\eta (X) =[X]$
by the definition of $\cH$-cycle classes.\\[4pt]

\B{rmk} Actually we got ``Gysin maps'' $i^!$ for
imbeddings  $i:X\into Y$ where $Y$ is just imbeddable in
$V$ smooth, by capping with the $\cH$-cycle class of $X$
in $V$. This operation will take a cycle of codimension
$p$ on $Y$ to a cycle on $X$ of codimension $p$ plus the
codimension of $Y$ in $V$.
\E{rmk}

Let denote $i_!:C_*(X)\to C_*(Y)$ the canonical map
induced by $i$. Since $Y$ is smooth we do have the
following equation
\B{equation} \eta i_! = i_{\diamond}\eta
\E{equation}
where $i_{\diamond}: H^*_X(Y,\cH^*)\to H^*(Y,\cH^*)$ is the
standard map. Let denote $i^*:  H^*(Y,\cH^*)\to
H^*(X,\cH^*)$. Let consider the intersection
product of cycles induced by the $\cH$-cohomology ring,
according with \S\ref{H-ring}.

\B{prop} The operation $i^!$ is functorial and compatible
with \'etale pull-backs. We have the self-intersection
property:
\B{equation}
i^!i_! (X) = X\bul X
\E{equation}
If $i:X\into Y$ is a smooth pair we then have
$$ i^! = i^*\eta \p{\cap} [X]$$
and $i^!$ is a ring homomorphism; there is a projection
formula $$i_!(x\bul i^!(y))=i_!(x)\bul y$$
for cycles $x$ and $y$ on $X$ and $Y$ respectively.
 \E{prop}
\B{proof} Compatibilities are easy to check. The
self-intersection property is obtained as follows:
\B{center}
\parbox{3in}{$i^!(i_!(x))=i_!(x)\p{\cap}\eta (X)=$\hfill
proj. form.\\$=x\p{\cap}i^*\eta (X)$}
\E{center}
whence, by taking $x=[X]$, we have
$$[X]\p{\cap}i^*\eta (X)=X\bul X $$
This last equation holds because of the
Scholium~\ref{prs}, giving us the following
$$\eta ([X]\p{\cap}i^*\eta (X))=\eta (X)\eta (X)$$
where $\eta (X)\eta (X)=\eta (X\bul X)$ (see Theorem~3)
and $\eta $ is an isomorphism.\\ The other equation is
given by the following commutative diagram
$$\begin{array}{ccc} H^{*}(Y,\cH^*)\tensor
H^{c}_X(Y,\cH^c\p{(c)}) &\longby{\cdot}  &
H^{*+c}_X(Y,\cH^{*+c}\p{(*+c)})\\ \p{i^*\tensor \cap
[Y]}\downarrow & &\downarrow \wr\ \p{\cap [Y]}\\
H^{*}(X,\cH^*)\tensor C_{d-c}(X;\La ) &
\longby{\cap}&  C_{d-c-*}(X;\La )
\end{array}
$$
which is obtained by the Scholium~\ref{prs} (cf. the
formula (\ref{rest}) ). The projection formula is obtained
from the projection formula w.r.t. the $\cH$-product (cf.
the Scholium~\ref{prs} and the Theorem~3).
 \E{proof}

By using the contravariant structure of
$\cH$-cohomologies we can construct `refined' Gysin maps
$f^!$ for algebraic cycles between imbeddable varieties.

\subsection{Grothendieck-Gillet axioms for Chern classes}
\label{H-chern}

The way to obtain a theory of Chern classes in
$\cH$-cohomologies and the corresponding Riemann-Roch
Theorems will be to show that the cohomology theory
$H^{*}_Z(X,\cH^*\p{(*)})$ and the homology theory
$C_{*,*}(-;\La\p{(*)})$ satisfy the list of axioms in
\cite[Definition 1.1 -- 1.2]{GIL}. We are going to
consider $X\in\cat{V}_k$ smooth over a perfect
field. We also assume that $\La$  (constant sheaf for the
Zariski topology) has finite weak global dimension (see
\cite[Definition 2.6.2]{KS}) in order to consider tensor
products $-\stackrel{L}{\otimes}_{\La}-$ in the derived
category.

\B{defi}\label{line} We let say that a natural
transformation  $$c\ell :\Pic (X)\otimes\La \to
H^1(X,\cH^1\p{(1)})\subset H^2(X,\p{1})$$ of contravariant
functors is a {\it cycle class map for line bundles}
if $c\ell$ localizes satisfiyng the properties
{\it (i) -- (ii)}\, stated in \S\ref{int-cyc}. Therefore
$c\ell$ is compatible via the local triviality property
\cite[1.5]{BO}, with the map obtained mapping a prime Weil
divisor $i: D\into X$ to the Poincar\'e dual of the direct
image under $i$ of the fundamental class $\eta_D$.
\E{defi}

\B{teor} Let $(H^*,H_*)$ be an appropriate Poincar\'e
duality on $\cat{V}_k$ for a perfect field $k$, with values
in a fixed category of $\La$-modules such that $H^*$
satisfies the homotopy property and there is a cycle class
map for line bundles. Then there is a theory of Chern
classes $$c_{p,i} : K_i^Z(X) \to
H^{p-i}_Z(X,\cH^p\p{(p)})$$ associated with any closed $Z$
in $X\in\cat{V}_k$ smooth.
\E{teor}
\B{proof} With the notations of \cite{GIL} we let
$\oplus\underline \Gamma^*\p{(p)}\df\oplus \cH^p\p{(p)}$
be the graded sheaf with the $\cH$-cup-product (according with
\cite[Definition 1.1]{GIL} and \S\ref{H-ring}) defining our
cohomology theory ring on the category $\cat{V}_k$.
We let define the homology as $$H_i(X,\Gamma \p{(j)})\df
\H^{-i}(X,{\cQ}^{j}_{\bul}\p{(j)})=C_{i,j}(X;\La\p{(j)})$$
which is covariant w.r.t. proper morphisms and a
presheaf for the \'etale topoloy by \cite[3.7]{BO}; the
compatibility \cite[1.2.(i)]{GIL} is ensured by the
compatibility \cite[1.2.2]{BO} and limits arguments (cf.
\S~4.3). The functorial long exact sequence of homology,
for a pair $i:Y\into X$, is obtained via the
hypercohomology long exact sequence with supports
$$\H^{-i}_Y(X,{\cQ}^{j}_{\bul}\p{(j)})\to
\H^{-i}(X,{\cQ}^{j}_{\bul}\p{(j)})\to
\H^{-i}(X-Y,{\cQ}^{j}_{\bul}\p{(j)})$$
since $\Gamma_Y{\cQ}^{j}_{\bul,X}\p{(j)}\cong
i_*{\cQ}^{j}_{\bul,Y}\p{(j)}$ and Lemma~4.6 (i.e.
\cite[1.2.(ii)]{GIL} holds). The cap product structure
is given by the $\cH$-cap-product and all the properties
required by \cite[1.2.(iii) -- (viii)]{GIL} are easily
seen by using the results of \S4 and \S5. The homotopy
property \cite[1.2.(ix)]{GIL} is our Lemma~3.4. Thus
we are left to show the following classical Dold-Thom
decomposition (see \cite[1.2.(x)--(xi)]{GIL}).
\E{proof}

\B{schol}\label{Edeco}{\rm  (Decomposition)} Let $\cE$ be a
locally free sheaf, $\rank \cE = n+1$, on $X$ smooth.  Let
$\pi: P\df {\bf P}(\cE) \to X$ be the corresponding
projective bundle. For $\cO_P(1)\in \Pic P$ let $$\xi \df
c\ell (\cO_P(1)) \in H^1(P,\cH^1\p{(1)})$$ we then have
 $$\oplus \pi^*( )\p{\cup}\xi^i: \bigoplus_{i=0}^{n}
H^{p-i}(X,\cH^{q-i}_X\p{(r-i)})\cong
H^p(P,\cH^q_P\p{(r)})$$
Furthermore: \B{equation}\label{trace}
\int_{P/X}^{}\xi^n =1_{\La}
\E{equation}
\E{schol}
\B{proof} The element $\xi\in H^1(P,\cH^1\p{(1)})$ defines
a map $\La \to \cH^1\p{(1)}[1]$ in the derived category
$\cat{D}(P_{Zar};\La)$ of complexes of sheaves of
$\La$-modules on the Zariski site. By cup-product we have
a map $\xi^i : \La [-i] \to \cH^i\p{(i)}$ hence a map
${\bf R}\xi^i :\La [-i] \to {\bf R}\pi_*\cH^i\p{(i)}$ in
$\cat{D}(X_{Zar};\La)$. On the other hand we have the
canonical map ${\bf R}\pi^{\sharp} :\cH^{q-i}\p{(r-i)} \to
{\bf R}\pi_*\cH^{q-i}\p{(r-i)}$ induced by contravariancy
(cf. \S 2.1). By taking
${\bf R}\pi^{\sharp} \stackrel{L}{\otimes} {\bf R}\xi^i$ we
obtain the maps  $$\cH^{q-i}\p{(r-i)}[-i]\to {\bf
R}\pi_*\cH^{q-i}\p{(r-i)} \stackrel{L}{\otimes}{\bf
R}\pi_*\cH^i\p{(i)}$$ Since the
$\cH$-cup-product yields products (cf.
\cite[Expos\'e \S2.1]{SGA5}) $${\bf
R}\pi_*\cH^{q-i}\p{(r-i)} \stackrel{L}{\otimes}{\bf
R}\pi_*\cH^i\p{(i)}\to {\bf R}\pi_*\cH^q\p{(r)}$$ by
composing we do get the maps  $\cH^{q-i}\p{(r-i)}[-i]\to
{\bf R}\pi_*\cH^q\p{(r)}$  whence the map  $$\gamma :
\bigoplus_{i=0}^{n}\cH^{q-i}\p{(r-i)}[-i]\to {\bf
R}\pi_*\cH^q\p{(r)}$$ in the derived category
$\cat{D}(X_{Zar};\La)$. Now that $\gamma$ is defined the
claimed decomposition will follows by proving that
$\gamma$ is a quasi-isomorphism because of the Leray
spectral sequence $$H^p(X,{\bf R}\pi_*\cH^q\p{(r)})\cong
H^p(P,\cH^q\p{(r)})$$
In order to show that $\gamma$ is a quasi-isomorphism we
are left to show the isomorphisms of groups
$$(\gamma^p)_x :
\bigoplus_{i=0}^{n}H^{p-i}({\rm
Spec}\cO_{X,x},\cH^{q-i}\p{(r-i)})\cong
(R^p\pi_*\cH^q\p{(r)})_x$$
for all $x\in X$ and $p\geq 0$. By continuity of the
arithmetic resolutions the stalks
$(R^p\pi_*\cH^q\p{(r)})_x$ are computed by
$H^{p}(\P^n_{\cO_{X,x}},\cH^{q}\p{(r)})$; we need the
following compatibility:
\B{lemma} Let $U\subset X$ be an
open Zariski neighborhood of $x$ on which $\cE$ is free.
Let $\xi\in  H^1_{\infty}(\P_U^n,\cH^1\p{(1)})$ be the
restriction of the tautological divisor, where $i:\infty
\cong \P_U^{n-1}\into \P_U^n$ is a hyperplane at infinity.
Then (with the notation of \S3.3)
$$j_{(n,n-1)}\pi^*_{n-1} = \pi^*_{n}\p{\cup}\xi$$
equality between maps from
$H^{p-1}(U,\cH^{q-1}\p{(r-1)})$ to
$H^p(\P_U^n,\cH^q\p{(r)})$.
 \E{lemma}
\B{proof}
Note that $\pi^*_{n-1}=i^*\pi^*_{n}$. The purity
isomorphism $$H^{p-1}(\infty ,\cH^{q-1}\p{(r-1)})\cong
H^p_{\infty}(\P_U^n,\cH^q\p{(r)})$$ is obtained as $\eta
(- \p{\cap} [\infty])$ and $\eta (\infty ) = \xi$
(because of the compabilities of the cycle
class $c\ell$) thus we have   \B{center}
\parbox{3in}{$\eta
(\pi^*_{n-1}\p{\cap}[\infty])=$\hfill\\  $=\eta
(i^*\pi^*_{n}\p{\cap}[\infty])=$\hfill Scholium~6.2\\
$=\eta (\pi^*_{n}\p{\cup}\eta
(\infty)\p{\cap}[\P^n_U])=$\hfill\\ $=\pi^*_{n}\p{\cup}\eta
(\infty) =$ \hfill\\ $=\pi^*_{n}\p{\cup}\xi$ \hfill as
elements in $H^p_{\infty}(\P_U^n,\cH^q\p{(r)})$}
\E{center}
and, by definition of $j_{(n,n-1)}$, the image of it
under $H^p_{\infty}(\P_U^n,\cH^q\p{(r)})\to
H^p(\P_U^n,\cH^q\p{(r)})$ yields the claimed equation.
\E{proof}
Thus: $(\gamma^p)_x$ is clearly an isomorphism by
reduction to open Zariski neighborhoods on which $\cE$ is
free, arguing as in \S3.3 via the Lemma above and
induction on the rank of $\cE$.\\  Let show the equation
(\ref{trace}). By the definition of the
$\cH$-Gysin map we have that $\pi_*: H^n(P,\cH^n\p{(n)})\to
H^0(X,\cH^0\p{(0)})$ is obtained (via the Leray spectral
sequence) by composition with $${\bf R}\pi_{\flat}:{\bf
R}\pi_*\cH^n\p{(n)}\to \cH^0\p{(0)}[-n]$$ in
$\cat{D}(X_{Zar};\La)$. Thus (see \S\ref{H-ring} ) we have
to prove that the composition of
$$\La [-n] \longby{{\bf R}\xi^n}{\bf
R}\pi_*\cH^n\p{(n)}\longby{{\bf
R}\pi_{\flat}}\cH^0\p{(0)}[-n]\cong\La [-n]$$ is the
identity. Arguing as above we are reduced to show the
equation (\ref{trace}) for $\pi:\P^n_{\cO_{X,x}}\to {\rm
Spec} \cO_{X,x}$. By the projection formula
$$\pi_*(\pi^*(1)\p{\cup}\xi^n)=\pi_*(\xi^n)$$
we are left to show that $\pi_*$ is the inverse of the
``decomposition'' isomorphism $$\pi^*( )\p{\cup}\xi^n :
H^0({\rm Spec} \cO_{X,x},\cH^0\p{(0)}) \to
H^n(\P^n_{\cO_{X,x}}, \cH^n\p{(n)})$$
By choosing a $k$-rational point of $\P^n_k$ we get a
proper section $\sigma$ of $\pi$. With the notation
above: $\sigma_* = j_{(n,0)}\pi^*_{0}$ and
by the Lemma we have $$\pi^*(1)\p{\cup}\xi^n
= \sigma_*(1)$$
By applying $\pi_*$ to the latter and taking the
image of it under the canonical augmentation $H^0({\rm
Spec} \cO_{X,x},\cH^0\p{(0)})\cong \La$ we do obtain the
claimed formula.
 \E{proof}
\B{rmk} After Grothendieck-Verdier, this
`decomposition argument' is quite standard. See
\cite[Expos\'e VII]{SGA5} for \'etale cohomology and
\cite[Theor.8.2]{GIL} or \cite{SH} for the $K$-theory.
 \E{rmk}
Let $$A(-) = \bigoplus_{p,q,r}^{}
H^p(-,\cH^q\p{(r)})$$
be the $\cH$-cohomology ring functor.
\B{cor} Let $\pi :\P(\cE)\to X$ be as above, $\rank \cE =
n+1$. Then
$$\pi^* : A(X) \to A(\P(\cE))$$
is an injective homomorphism of unitary $\La$-algebras
and the elements $$1,\xi , \ldots , \xi^n$$ generate freely
$A(\P(\cE))$ as $A(X)$-module. Furthermore
$$\pi_* : A(\P(\cE)) \to A(X)$$
is a surjective homomorphism of $A(X)$-modules (having
degree $-n$).
\E{cor}
\B{proof} The statement is clear after \S\ref{H-ring} and
the Scholium above. For example, $\pi_*(\xi^i)=0$ for
$i=0,\ldots ,n-1$ but $\pi_*(\xi^n)=1$ by (\ref{trace})
whence the linear independence of $1,\xi , \ldots , \xi^n$
can be seen as follows: let suppose that $$\pi^*(x_0)+
\cdots +\pi^*(x_n)\p{\cup}\xi^n =0$$ then by applying $\pi_*$
and the projection formula we get $x_n=0$ thus
$$\pi^*(x_0)\p{\cup}\xi+ \cdots +\pi^*(x_{n-1})\p{\cup}\xi^n =0$$
and the same argument gives $x_{n-1}=0$ and so on.
Again: $\pi_*$ is a surjection because of
$$\pi_*(\pi^*(\dag )\p{\cup}\xi^n)=\pi_*(\xi^n)\p{\cup}\dag
=\dag$$
 \E{proof}

\B{rmk} By the prescription of \cite{GC} we therefore
obtain Chern classes $c_p : K_0(X) \to H^p(X,\cH^p\p{(p)})$
satisfying the equation
$$\xi^{n}+\pi^*c_{1}(E)\xi^{n-1}+\cdots +\pi^*c_n(E)=0$$
for $E$ a vector bundle of rank $n$. By \cite{GIL} we
have that $c_p$ is just $c_{p,0}$.
 \E{rmk}

\subsection{Variation on the invariance theme}
\label{finite}

Let consider a sophisticated Poincar\'e duality theory
$(H^*,H_*)$ satisfying the point axiom. Let consider
$f:X\to Y$ a proper dominant morphism between connected
smooth schemes in $\cat{V}_k$. If dim$X$ = dim$Y$ then
$K(X)$ is a finite field extension of $K(Y)$; let ${\rm
deg} f=[K(X):K(Y)]$ be its degree. Following the proof of
(\ref{prfor}) we have that the composition
 $$H^{*}(Y,\cdot) \by{f^{\star}}
H^{*}(X,\cdot) \by{f_{\star}}
H^{*}(Y,\cdot)$$ is the multiplication by ${\rm deg} f$,
as a consequence of the projection formula and our
assumption that $f_!(\eta_X)={\rm deg} f\cdot\eta_Y$
(cf. \S\ref{inter1}). Thus:
\B{prop}\label{mult} The composition of
$$H^p_Z(Y,\cH^{q}\p{(r)}) \by{f^{*}}
H^p_{f^{-1}(Z)}(X,\cH^{q}\p{(r)}) \by{f_{*}}
H^p_Z(Y,\cH^{q}\p{(r)}) $$
is the multiplication by ${\rm deg} f$.
\E{prop}
\B{proof} Let $d=$dim$X=$dim$Y$. Then the projection
formula (\ref{hpr}) looks
$$\B{array}{c}
\hspace{10pt} f_*{\cQ}^d_{\bul}\p{(d)}\tensor
f_*\cH^{q}\p{(r)}\\ \p{id\tensor f^{\sharp}}\nearrow \
\hspace{40pt} \searrow \p{f_*\cap_{{\cH}}}  \\
f_*{\cQ}^d_{\bul}\p{(d)}\tensor
\cH^{q}\p{(r)}\hspace{60pt}
f_*{\cQ}^{d-q}_{\bul}\p{(d-r)}\\  \p{f_{\sharp}\tensor
id}\downarrow \hspace{105pt}\downarrow \p{f_{\sharp}} \\
{\cQ}^d_{\bul}\p{(d)}\tensor \cH^{q}\p{(r)}
\hspace{10pt} \longby{\cap_{{\cH}}}
\hspace{10pt}{\cQ}^{d-q}_{\bul}\p{(d-r)}
\E{array}$$
By the dimension axiom the complex
${\cQ}^d_{\bul}\p{(d)}$ is concentrated in degree $d$ and
its hypercohomology $C_{d,d}(X,\La\p{(d)})$ has a natural
global section $[X]$ corresponding to the fundamental
class $\eta_X\in H_{2d}(X,\p{(d)})$. The same holds on
$Y$ and $[X]\leadsto {\rm deg}f [Y]$ under $f_*$. Thus by
taking cohomology with supports we have the  result.
 \E{proof}

\B{rmk} The same argument applies to the $K$-theory by
using the projection formula in (\ref{kpr}).
\E{rmk}

\B{lemma}\label{birsplit} Let $f:X'\to X$ be a proper
birational morphism between smooth varieties; let $i:
Z\into X$ and $i':Z'=f^{-1}(Z)\into X'$  be closed
subschemes such that $f: X'-Z' \cong X-Z$. Then we have
splitting short exact sequences $$0\to
H^p_Z(X,\cH^q\p{(r)}) \by{u} H^p(X,\cH^q\p{(r)})\oplus
H^p_{Z'}(X',\cH^q\p{(r)}) \by{v}
 H^p(X',\cH^q\p{(r)})\to 0$$
where:
$$ u=\left( \begin{array}{c}
  {i_{\diamond}}\\{f^*}
\end{array}     \right) $$
and
$$v= (f^*,-i'_{\diamond})$$
The left splitting of $u$ is given by $u':(0,f_*)$.
 \E{lemma}
\B{proof} Let consider the following maps of long exact
sequences $$ \begin{array}{ccccccccc}
\cdots & \to & H^p_{Z'}(X',\cH^q\p{(r)}) &
{\by{i'_{\diamond}}}& H^p(X',\cH^q\p{(r)}) & \to &
H^{p}(X'-Z',\cH^q\p{(r)}) & \to & \cdots \\
 & &{\p{f_*}\downarrow\uparrow\p{f^*}} &
&{\downarrow\uparrow}  & & {\downarrow\wr} & &\\
{\cdots}&{\to}& H^p_Z(X,\cH^q\p{(r)}) &{\by{i_{\diamond}}}&
H^p(X,\cH^q\p{(r)}) & {\to}&
H^{p}(X-Z,\cH^q\p{(r)}) &{\to}&{\cdots}
\end{array} $$
Since deg$f=1$ by the Proposition above $f_*f^*=1$;
thus the corresponding Mayer-Vietoris exact sequence
splits (because the boundary is zero) in short exact
sequences as claimed.
\E{proof}

\B{schol} Let $X$, $X'$, $Z$ and $Z'$ be as above and
pure dimensional. We have isomorphisms ($d=$dim$X$)
$$C_{d-p}(X;\La)\oplus C_{d-p}(Z';\La) \by{\simeq}
C_{d-p}(Z;\La)\oplus
C_{d-p}(X';\La)$$  given by the matrix
$$\left(
\B{array}{cc} 0 & f_! \\ f^! & -i'_!
\E{array}\right) $$ and, for $Z$ and $Z'$ smooth of
codimension $c$ and $c'$: $$H^{p}(X,\cH^{q}\p{(r)})\oplus
H^{p-c'}(Z',\cH^{q-c'}\p{(r-c')}) \by{\simeq}
H^{p-c}(Z,\cH^{q-c}\p{(r-c)})\oplus
H^{p}(X',\cH^{q}\p{(r)})$$  given by the matrix
$$\left(\B{array}{cc} 0 & f_* \\ f^* & -j_{\pp{(X',Z')}}
\E{array}\right) $$
where $j_{\pp{(X',Z')}}$ is the Gysin map in
$\cH$-cohomology (cf. Scholium~\ref{Gysin},
\S\ref{H-Gys}). \E{schol}
\B{proof}
By the Lemma~\ref{birsplit} and purity.
\E{proof}

\subsection{Blowing-up}

Let $(H^*,H_*)$ be an appropriate Poincar\'e
duality such that $H^*$ satisfies the homotopy property and
there is a cycle class map for line bundles (cf.
\S\ref{H-chern}).\\  Let $f:X'\to X$ be the blow up of a
smooth subvariety $Z$ of codimension $c\geq 2$ in a smooth
variety $X$ of dimension $d$. Thus the exceptional divisor
is the projective bundle over $Z$ given by $\P (\cN)$
where $\cN$ is the normal sheaf, locally free  of rank $c$.

\B{prop}\label{blow-up} For the blow-up $X'$ of $X$ along
$Z$ as above we have the following canonical formulas
$$H^p(X',\cH^q\p{(r)}) \cong H^p(X,\cH^q\p{(r)})\oplus
\bigoplus_{i=0}^{c-2}
H^{p-1-i}(Z,\cH^{q-1-i}\p{(r-1-i)})$$
and in particular
$$C_{n}(X';\La) \cong C_{n}(X;\La) \oplus
\bigoplus_{i=0}^{c-2} C_{n-c+1+i}(Z;\La)$$
\E{prop}
\B{proof} Since $f:Z'\to Z$ is a proper morphism,
between smooth varieties, having relative dimension $1-c$
 we have a push-forward (see \S\ref{H-Gys}) $$f_*:
H^{p-1}(Z',\cH^{q-1} \p{(r-1)}) \to
H^{p-c}(Z,\cH^{q-c}\p{(r-c)})$$ Because
of purity $f^*: H^p_Z(X,\cH^q\p{(r)}) \into
H^p_{Z'}(X',\cH^q\p{(r)})$ induces a map
$$f^!: H^{p-c}(Z,\cH^{q-c}\p{(r-c)}) \into
 H^{p-1}(Z',\cH^{q-1}\p{(r-1)})$$ as well. By the
Lemma~\ref{birsplit} we have a splitting exact sequence
$$0\to H^{p-c}(Z,\cH^{q-c}\p{(r-c)}) \by{u}
H^p(X,\cH^q\p{(r)})\oplus H^{p-1}(Z',\cH^{q-1}\p{(r-1)})
\by{v} H^p(X',\cH^q\p{(r)})\to 0$$
with left splitting $u'= (0,f_*)$. Let consider the
projector $\pi =uu'$; we then have $\pi u=u$, $v\pi =0$
thus $v$ restricts to an isomorphism
$$v: {\rm ker} \pi \cong  H^p(X',\cH^q\p{(r)})$$
Now $\pi (x, z') = u (f_*(z')) = (j_{\pp{(X,Z)}}(f_*(z')),
f^!f_*(z'))=0$  if and only if $f_*(z')=0$ (because $f^!$
is injective). Thus we have
$$ {\rm ker} \pi = H^p(X,\cH^q\p{(r)})\oplus {\rm ker}
f_*$$
Since $Z' =\P (\cN)$ and $f_{\mid Z'}$ is the standard
projection then, by the Dold-Thom decomposition (see
Scholium~\ref{Edeco}), we have an exact sequence
$$0\to\bigoplus_{i=0}^{c-2}
H^{p-1-i}(Z,\cH^{q-1-i}\p{(r-1-i)}) \by{\xi^i f^*}
H^{p-1}(Z',\cH^{q-1} \p{(r-1)}) \by{f_*}
H^{p-c}(Z,\cH^{q-c} \p{(r-c)})\to 0 $$  where $\xi$ is the
tautological divisor, hence:
 $${\rm ker} \pi \cong H^p(X,\cH^q\p{(r)})
\oplus\bigoplus_{i=0}^{c-2}
H^{p-1-i}(Z,\cH^{q-1-i}\p{(r-1-i)})$$
and the claimed isomorphisms are easily obtained.
 \E{proof}

\B{rmk} For the $\cK$-cohomology the same proof
applies yielding the formula $$H^p(X',\cK_q) \cong
H^p(X,\cK_q)\oplus  \bigoplus_{i=0}^{c-2}
H^{p-1-i}(Z,\cK_{q-1-i})$$
\E{rmk}

\appendix

\section{Examples and comments}\label{app}

We will give a draft for testing the common cohomologies
with respect to our setting (cf. \S\ref{H-cap},
\S\ref{H-cup}).

\subsection{Grothendieck-Verdier duality}

Let $X_{fine}$ be a site finer than the Zariski site
$X_{Zar}$ for $X\in\cat{V}_k$ e.g. \'etale or analytic
sites. Let $\La$ be a commutative ring with $1$ (of
finite weak global dimension) and let $\cat{D}(X;\La)$
denote the derived category of the abelian category of
complexes of $\La$-modules in the corresponding
Grothendieck topos on $X_{fine}$. We let $F^{\bul}\leadsto
F^{\bul}\p{(r)}$ denote a ``twist \`a la Tate'' functor on
$\cat{D}(X;\La)$ commuting with direct and inverse image
functors. For a closed imbedding $i:Z\into X$ let
$j:U=X-Z\into X$ be the corresponding Zariski open
immersion; let consider the six standard operations
$$i^*\dashv i_* \dashv i^! \ \ \ \  j_!\dashv j^* \dashv
Rj_*$$ where $i^!$ is the ``sheaf of sections on $Z$''
functor  and $j_!$ is the ``extension by zero'' functor.
Let consider twisted objects $\La^{\bul}_k\p{(r)}\in
\cat{D}(k;\La)$ with a canonical augmentation (usually a
quasi-isomorphism) $\La\to\La^{\bul}_k\p{(0)}$ such that
$\oplus_{r}^{}\La^{\bul}_k\p{(r)}$ is a graded
$\La$-algebra via the canonical maps
$$m:\La^{\bul}_k\p{(r)}\stackrel{L}{\otimes}
\La^{\bul}_k\p{(s)}\to \La^{\bul}_k\p{(r+s)}$$ We let
extends $\oplus_{r}^{}\La^{\bul}\p{(r)}$  to a
$\La$-algebra over the big site $(\cat{V}_k)_{fine}$ by
pulling back along the structural morphisms. Thus we have
a contravariant functor on $\cat{V}_k^2$  $$(X,Z)\leadsto
Hom_X(\La , i_*i^!\La^{\bul}_X\p{(r)}[q])\df
H_Z^q(X,\p{r})$$ where the $Hom$ is taken in
$\cat{D}(X;\La)$. By the general non-sense of triangulated
categories this is a cohomology theory with cup-products:
the long exact sequence of cohomology with supports is
given by the triangles (since the $Hom$ of triangulated
categories takes triangles to long exact sequences)
$$j_!j^*\La \to \La \to i_*i^*\La \by{+1}$$ and the
cup-product of $a: \La \to \La^{\bul}\p{(r)}[q]$ and $b:
\La \to \La^{\bul}\p{(r')}[q']$ is given by tensoring
$\La^{\bul}\p{(r)}[q]\otimes b: \La^{\bul}\p{(r)}[q]\to
\La^{\bul}\p{(r)}\stackrel{L}
{\otimes}\La^{\bul}\p{(r')}[q+q']$ and composing $a$ with
$$\La^{\bul}\p{(r)}[q]\to
\La^{\bul}\p{(r)}\stackrel{L}{\otimes}\La^{\bul}\p{(r')}
[q+q'] \by{m[q+q']} \La^{\bul}\p{(r+r')}[q+q']$$ in the
derived category. (Note: this is the natural $Ext$
pairing). Since the composition of maps is associative it
is clear that the pairing above is compatible with the
long exact sequences of cohomology with supports and
furthermore the following square commutes $$\B{array}{ccc}
Hom_X(j_!j^*\La ,\La^{\bul}\p{(r)}[q])\otimes
H^{q'}(X,\p{r'})&\to &Hom_X(j_!j^*\La ,
\La^{\bul}\p{(r+r')}[q+q'])\\  \p{\theta\otimes
j^*}\downarrow& &\downarrow\p{\theta}\\ Hom_U(\La
,\La^{\bul}\p{(r)}[q])\otimes H^{q'}(U,\p{r'})&\to
&Hom_U(\La , \La^{\bul}\p{(r+r')}[q+q']) \E{array} $$
where $\theta$ is the canonical isomorphism induced by
$j_!\dashv j^*$ and the bottom arrow is the cup-product
pairing on $U$: thus the axioms $\forall 1-\forall 3$ (see
\S5) are  satisfied. (Note: to check commutativity of the
square above we simply need that $j^*(b)\p{\circ}\theta
(a) = \theta (b\p{\circ}a)$ which is a consequence of
$j_!\dashv j^*$).\\ Let assume the existence of a ``global
duality'' $$f_!\dashv f^!$$ where: $f_! = Rf_*$ for $f$
proper,  $f^!=f^*$ for $f$ \'etale (or at least a Zariski
open imbedding) and $f^! =f^*\p{(d)}[2d]$ for $f$ smooth of
relative dimension $d$. Thus we obtain a homology theory
$$ X \leadsto Hom_X(\La ,
\pi^!\La^{\bul}_k\p{(-r)}[-q])\df H_q(X,\p{r})$$ where
$\pi : X \to k$ is the structural morphism. Furthermore
the counit $f_!f^! \to 1$ yields a pairing
$$f^*(\La^{\bul}\p{(r)})\stackrel{L}{\otimes}
f^!(\La^{\bul}\p{(r')})\to
f^!(\La^{\bul}\p{(r)}\stackrel{L}{\otimes}\La^{\bul}\p{(r')})
\by{f^!(m)}f^!(\La^{\bul}\p{(r+r')})$$ whence applied to
$\pi$ and the previous procedure give us a sophisticated
cap-product and the axioms A1--A4 (see \S4) are clearly
satisfied.\\ For $\pi :X \to k$ smooth of pure dimension
$d$ Poincar\'e duality will mean to be given by natural
quasi-isomorphisms: $$\eta :\pi^*\La^{\bul}_k\p{(0)}
\by{\simeq} \pi^!\La^{\bul}_k \p{(-d)}[-2d] \mbox{\ \ \
``fundamental class''}$$ (yielding a homology class by
composing with the augmentation) and
$$\pi^*\La^{\bul}_k\p{(r)}[q] \by{\simeq}
\pi^!\La^{\bul}_k \p{(r-d)}[q-2d]$$ Thus the compatibility
of cap and cup products (see \S6) is ensured by
construction. Concerning the dimension axiom (see A5--A6
in \S4.4) just compare with the Lemma 2.1.2 of
\cite{BO}.\\  This framework applies to: {\it (i)}\, the
analytic site where $\La =\Z$ and  $\La^{\bul}_k\p{(r)}\df
\Z\p{(r)}$ by mean of the Tate twist in Hodge theory (or
$\La \df {\bf Q}, {\bf R}$ or complex-De Rham cohomology
where $\La =\La^{\bul}_k\p{(r)} \df \C$) after Verdier
Expos\'e \cite{VE} (cf. \cite{KS}), {\it (ii)}\, the
\'etale site where $\La \df \Z/\ell^{\nu}$ and
$\La^{\bul}_k\p{(r)}\df \Z/\ell^{\nu}(r)$= \'etale sheaf
$\mu_{\ell^{\nu}}$ of $\ell^{\nu}$-th roots of unity
($\ell$ prime to char$k$) with Tate twist $r$, after
Deligne Expos\'e XVIII in \cite{SGA4}.

\subsection{Algebraic De Rham cohomology}

Let $X$ be a smooth algebraic scheme over $k=\C$. On the
analytic manifold $\aX$ the constant sheaf $\La =\C$ is
quasi-isomorphic to the holomorphic De Rham complex
$\Omega^{\bul}_{an}$; by the Grothendieck comparison
theorem the {\it algebraic}\, De Rham complex
$\Omega^{\bul}_{X/k}$ computes the complex cohomology as
well. In \cite{HA1},\cite{HA2} has been developed a
general theory of algebraic De Rham complexes for
imbeddable varieties over $k$ of characteristic zero. This
theory yields cohomology groups $H^*_{DR}(X)$ and homology
groups $H_*^{DR}(X)$ satisfying the Bloch-Ogus axioms (see
\cite[2.2]{BO}) e.g. $H_*^{DR}(X)\df
\H^{2d-*}_X(Y,\Omega^{\bul}_{Y/k})$ for $X\into Y$ and
$Y$ smooth of dimension $d$. Furthermore, the exterior
algebra structure on the De Rham complex grant us of
compatible cup products and sophisticated cap products as
explained in \cite[II.7.4]{HA2}. If moreover $k$ is
algebraically closed then $H^0_{DR}(X)=k$ for $X$
connected (see  \cite[II.7.1]{HA2}); since $H_q^{DR}(X)=0$
if $q>2d$ (see  \cite[II.7.2]{HA2}) the dimension axiom
(see \S4.4) is satisfied. The homotopy axiom is ensured by
the Poincar\'e Lemma (see \cite[II.7.1]{HA2}). Thus
algebraic De Rham cohomology and homology is appropriate
(in the sense of \S6.1) where $\La = k$ is algebraically
closed of characteristic zero. What about positive
characteristics?

\subsection{Deligne-Beilinson cohomology}

By considering the augmented and truncated complexes
$\cD\p{(r)}\df\Z \p{(r)} \to \Omega^{\bul<r}_{an}$ on the
complex compact manifold $X_{an}$ we let define the
Deligne-Beilinson cohomology $H^q_{\cD}(X;\Z \p{(r)})$ to
be the hypercohomology groups $\H^q(X_{an},\cD\p{(r)})$.
This definition can be `algebraically' extended: {\it
(i)}\, to smooth open varieties by taking
compactifications with normal crossing divisors at
infinity by mean of differential forms with logarithmic
poles at infinity, and {\it (ii)}\, to singular varieties
by smooth hypercoverings, according with Deligne-Hodge
theory. We refer to \cite{HV} for a detailed exposition of
the smooth case and \cite{GID}, \cite{JAD}, \cite{LED} for
the singular case: indeed, by using currents and
$\cC^{\infty}-$chains, we do have Deligne homology groups
$H_*^{\cD}(X;\Z \p{(\cdot)})$ and therefore a Poincar\'e
duality theory in the sense of Bloch-Ogus \cite[Theorem
1.19]{JAD}, \cite{GID}, \cite{BE}. Because of \cite[Lemma
1.20]{JAD} Deligne homology satisfies the dimension axiom.
Because of \cite[1.2 and \S3]{HV} we have products
$\p{\cup}:  \cD\p{(r)}\stackrel{L}{\otimes}\cD\p{(r')}\to
\cD\p{(r+r')}$ whence $H^*_{\cD}(-;\Z \p{(\cdot)})$ is
multiplicative (since the cup-product is obtained in a
standard way via the $Ext$ pairing and it is therefore
compatible with the long exact sequence of cohomology with
supports). By making the Exercise 1.8.6 in \cite{BE} one
obtains a sophisticated Poincar\'e duality which is
appropriate for algebraic cycles. (Note: $H^p(X,
\cH^p_{\cD}(\Z \p{(p)}))\cong CH^p(X)$ for $X$ smooth by
\cite{GID}). The homotopy property is obtained by those of
De Rham theory and integral cohomology via the long exact
sequences \cite[2.10]{HV} (cf. \cite[8.5]{HV}). (Note: we
are regarding the homotopy property as a property of mixed
Hodge structures).

\subsection{Cycle class map for line bundles}

Because of the Definition~7.2 (cf. the assumptions in
\S6.3) we let explain briefly some examples. For the sake
of exposition we first make the following remark.

\subsubsection{Cycle classes are locally trivial}

Roughly speaking, another way to understand
globally what such compatibility really mean is the
following. Let $X$ be smooth over a field and let assume
that our appropriate cohomology is coming with a
functorial map $c\ell :\Pic (X)\to H^2(X,\p{1})$ as usual
it is. Then the image of $c\ell$ is contained in the
subgroup of locally trivial elements i.e. the kernel of
$H^2(X,\p{1})\to H^0(X,\cH^2\p{(1)})$, simply because the
$\Pic$ of a local ring is zero. By the coniveau spectral
sequence and the point axiom we have that
$H^1(X,\cH^1\p{(1)})$ is identified with the subgroup of
locally trivial elements in $H^2(X,\p{1})$. By
construction of the $\cH$-cap-products we have the
following commutative square ($d=$dim$X$)
\B{displaymath}\B{array}{ccc}
H^{1}(X,\cH^1\p{(1)}) & \into &
H^2(X,\p{1}) \\  \p{\cap [X]}\ \downarrow& &
\downarrow\ \p{\cap \eta_X} \\ C_{d-1}(X;\La) &\by{\zeta} &
H_{2d-2}(X,\p{d-1})  \E{array}\E{displaymath}
where $\zeta$ is the cycle map in homology as defined in
\cite[\S7]{BO}. Thus the commutativity of
\B{displaymath}\B{array}{ccc}
\Pic (X) & \by{c\ell} &
H^2(X,\p{1}) \\  \p{w}\ \downarrow& &
\downarrow\ \p{\cap \eta_X} \\ CH_{d-1}(X) &\to
& H_{2d-2}(X,\p{d-1})  \E{array}\E{displaymath}
turns out to be an equivalent formulation of our
requirement (where $w$ is the associated Weil divisor
mapping.) To obtain the corrresponding local formulation
one has to make use of the supports of a given line bundle.

\subsubsection{Classical integral cohomology}

Let $\omega : \aX \to X_{Zar}$ be the canonical
continuous map of sites, where $X$ is a reduced
irreducible algebraic scheme over $\C$ and let identify
$\cH^{q}\p{(r)}$ by the Zariski sheaf
$R^{q}\omega_{\ast}\Z\p{(r)}$. From the exponential
sequence on $\aX$ we get a map $\omega_*
\cO_{\sa}^{\ast}\to \cH^{1}\p{(1)}$. Since any regular
mapping is analytic, $\cO_{X}^{\ast}$ is a subsheaf of
$\omega_{\ast}\cO_{\sa}^{\ast}$ whence the cycle map
$c\ell : \Pic (X) \to H^1(X,\cH^{1}\p{(1)})$ is obtained
by applying $H^1_{Zar}(X,-)$ to the boundary map
$\cO_{X}^{\ast}\to \cH^{1}\p{(1)}$. (Note: because of the
identifications
$H^0(X,\cH^{1}\p{(1)})=H^1(\aX,\Z)=Hom(H_1(\aX),\Z)$ this
boundary map takes the continuous mapping $f:\aX\to\C^*$
to the mapping that to a given loop $\phi:S^1\to\aX$
associates the Brouwer degree of the composition
$S^1\by{\phi}\aX\by{f}\C^*\by{n}S^1$ where $n(z)\df
\frac{z}{\mid z\mid}$ is the canonical homotopical
retraction. For $X$ affine hence $\aX$ a Stein space this
is a surjection as a map from the holomorphic functions on
$\aX$. The non-vanishing of it, for $f$ holomorphic
non-constant, computes the obstruction of having an
holomorphic logarithm according with the celebrated
Riemann representation theorem cf. \cite[V.\S2.4-5]{GRE}.)
To see that this cycle map $c\ell$ is compatible with the
classical $\Pic (\aX) \to H^2(\aX,\Z)$ we refer to
\cite[Proposition 1]{BS}. Because of the functoriality of
the exponential sequence $c\ell$ is a natural
transformation of contravariant functors meeting our
hypothesis {\it (i)}\, of \S6.3 if taking $H^1_D(X,-)$. To
check the assumption {\it (ii)}\, ibid. we let remark that
the map $\cO_{X}^{\ast}\to \cH^{1}\p{(1)}$ yields maps on
the stalks $\cO_{X,x}^{\ast}\to H^1(x)$ at $x\in X$ and in
particular at the generic point give us a map  $K(X)^*\to
H^1(K(X))$. By capping with the fundamental class we then
have a map $K(X)^*\to H_{2d-1}(K(X),\p{d-1})$
($d=$dim$X$); in the same way the rank mapping give us the
canonical map  $K_0 \to H_{2d}(-,\p{d})$ which is a local
isomomorphism. We then have a map between augmented
``Gersten complexes''  \B{displaymath}\B{array}{ccccc}
\cO_{X}^{\ast}&\to&K(X)^*&\by{div}&\coprod_{x\in
X^1}^{}i_x\Z\\ \downarrow&&\downarrow&&\mid\mid \\
\cH^{1}\p{(1)}&\to&H_{2d-1}(K(X),\p{d-1})&\by{d}&\coprod_{x\in
X^1}^{}i_x\Z
 \E{array}\E{displaymath}
where the bottom augmentation is given by restriction to
the generic point and cap product with the fundamental
class. (Warning: to check that this is actually a map of
complexes we should have to show that the divisor map $div$
has the claimed factorisation via $d$, the differential in
the niveau spectral sequence. This is ensured by the
Tautology 7.2 and the Theorem 7.3 in \cite{BO}; note
that the non-singular case grant us for the general case by
resolution of singularities.) Thus, by applying
$H^1_D(X,-)$ to these augmentations, we obtain the
commutative square in {\it (ii)}\, of \S6.3 as desired.

\subsubsection{The \'etale theory}

A similar framework as above applies to the \'etale
cohomology of  the sheaf $\mu_{\ell^{\nu}}$ by using the
Kummer sequence on the \'etale site ($\ell\neq$char$k$)
where the cycle map $\Pic (X)\otimes\Z/\ell^{\nu}\to
H^1(X,\cH^1(\mu_{\ell^{\nu}}))$ is now induced by the
identification $\cO_{X}^{\ast}\otimes Z/\ell^{\nu}\cong
\cH^1(\mu_{\ell^{\nu}})$ obtained by the Hilbert's
theorem 90 and the well known fact that the Picard group
of a local ring is zero; it is an analogous exercise to
check the aimed properties (by using \cite[Theorem
7.7]{BO}).

\subsubsection{Algebraic De Rham cohomology}

In the case of algebraic De Rham cohomology a nice
exposition of these properties has been written by
Hartshorne in \cite[II.7.6-7]{HA2}.

\subsubsection{Deligne-Beilinson cohomology}

For the Deligne-Beilinson cohomology $H^*_{\cD}(X;\Z
\p{(\cdot)})$ we dispose of a canonical map $$\rho
:H^0(X,\cO^*_X) \to  H^1_{\cD}(X;\Z \p{(1)})$$ (see
\cite[1.4.ii) and 2.12.iii)]{HV}, \cite[(1.1)]{LED}). By
sheafifying for the Zariski topology we do get $\rho :
\cO^*_X\to  \cH^1_{\cD}(\Z \p{(1)})$ which is an
isomorphism on $X$ smooth thus: $c\ell : \Pic (X) \cong
H^1(X,\cH^1_{\cD}(\Z \p{(1)}))$. Moreover: for any $X$  we
have that the kernel of the divisor map $div$ is given by
$H_{2d-1}^{\cD}(X;\Z\p{(d-1)})$ by \cite[Lemma 3.1]{JAD}.
Thus, by arguing as in the case of integral cohomology,
we easily obtain the aimed properties. (Note: the diagram
(3.1.1) of \cite{JAD} contains implicitly our
requirement.) See also \cite[1.9.2]{BE}.

\subsection{Singular varieties and intersection cohomology}

In order to obtain an ``intersection theory'' for singular
varieties, after \S5.5 and Corollary~6.6, one is tempted to
make use of $\cH$-cohomology rings. For complex
algebraic varieties Deligne-Beilinson cohomology do most
of the job but it is an open question whether the
$\cH$-cohomology ring functor is covariant for proper
morphisms. Another problem: the existence of
$\cH$-cycle classes on singular varieties. (Note: a
necessary condition is local triviality.) By weakening to
intersections modulo algebraic equivalence the same
problems occur with singular cohomology. Furthermore it
is hopeless to expect birational invariance of
$H^0(X,\cH^*)$ e.g. because of the following
counterexample: let $S\subset \P^3$ be the complex surface
defined by the homogeneous equation
$w(x^3-y^2z)+f(x,y,z)=0$ where $f$ is a general
homogeneous polynomial of degree $4$, then $S$  is a
rational quartic with a triple point; by the
computations of \cite{BSM} $H^1(S,\cH^1(\Z))\neq
H^2(S_{an} ,\Z)$ whence $H^0(S,\cH^2(\Z))\neq 0$  but
$H^0(S',\cH^2(\Z))= 0$ for any resolution of singularities
$f:S'\to S$ since $S'$ is rational; moreover
$H^0(S,\cH^2(\Z))$ is free thus
$H^0(S,\cH^2(\Z))\otimes\Z/n \subset
H^0(S,\cH^2(\mu_n))\neq 0$ which implies that $S$ has
non-zero torsion Brauer group. Similarly:
$\cH$-cohomologies don't satisfy the homotopy axiom in
the singular case.\\ It might be interesting the
approach to the ``intersection theory'' for singular
varieties by considering intersection cohomologies. For
$\aX$ normal we dispose of nice cohomology groups
$IH^*(\aX)$(= the intersection cohomology with middle
perverity, see \cite{BBD} and \cite{GOM}) having an
``intersection pairing'' and self duality (rationally).
By sheafifying for the Zariski topology we do have a sheaf
$\cI\cH^*$ and the corresponding Cousin complex associated
with the filtration by codimension (see \cite{HA1}).
(Note: this is the `homological' Gersten complex if local
purity holds.) So far as I'm concerned it is natural to
ask if $\cI\cH^*$ is Cohen-Macaulay i.e. is the Cousin
complex a flasque resolution? A positive answer would be a
great improvement in matter.

\subsection{Motivic cohomology}

Following Grothendieck philosophy of `motives' we are
suspecting the existence of an appropriate Poincar\'e
duality theory $H^*_{\cM}(-,\Z \p{(\cdot)})$ satisfying
the homotopy axiom for smooth algebraic varieties which is
`universal' at least in the following weak meaning (cf.
\cite[\S0.2]{BMS}). Let $\cat{P}\cat{V}_k$ be the
category whose objects are appropriate dualities (may be
with values in an abelian tensor category not just graded
abelian groups, see \cite[\S6]{JA}) and morphisms are
natural transformations compatible with products. For
example we have a canonical morphism from the
Deligne-Beilinson theory to the classical Borel-Moore
theory. We are expecting that any object
$H^*\in\cat{P}\cat{V}_k$ is a receptor for a map from
`motivic cohomology'  $H^*_{\cM}(-,\Z \p{(\cdot)})$. This
is to say that `motivic cohomology' is to some extent
`initial' in  $\cat{P}\cat{V}_k$. A reasonable candidate
for $H^q_{\cM}(-,{\bf Q} \p{(r)})\df H^q_{\cM}(-,\Z
\p{(r)})\otimes {\bf Q}$ is
$gr_{\gamma}^rK_{2r-q}(X)\otimes {\bf Q}$ after
Beilinson-Soul\'e, see \cite[Theorem 9]{SOU}. (Note: our
dimension axiom is the Beilinson-Soul\'e conjecture
\cite[2.2.2]{BE}). We would remark that if such
`motivic cohomology' $H^*_{\cM}(-,\Z \p{(\cdot)})$ exists
then the $\cH_{\cM}$-cohomology ring will be the Chow
ring since rational equivalence is the finer adequate
equivalence relation, see \cite[Prop. 8]{SAM}. Finally,
if we let restrict our attention to `representable'
dualities i.e. given by the homology of complexes, we
would remark that our cohomologies would have the descent
property of \cite{TH} since the axioms \cite[\S1.1]{BO}
ensured the Mayer-Vietoris property.

\B{thebibliography}{}

\bibitem[SGA 4]{SGA4}{\sc A.~Grothendieck~}et~al.:
Th\'eorie des topos et cohomologie \'etale des sch\'emas
(1963--64) {\it Springer} LNM {\bf 269 270 305},
Heidelberg, 1972-73.

\bibitem[SGA 5]{SGA5}{\sc A.~Grothendieck~}et~al.:
Cohomologie $\ell$-adique et fonctions L (1965--66)
{\it Springer} LNM {\bf 589}, Heidelberg, 1977.

\bibitem{BV1}{\sc L.Barbieri-Viale}: Des invariants
birationnels associ\'es aux th\'eories
cohomologiques, {\it C.R.Acad.Sci.Paris} {\bf 315}
S\'erie I (1992), 1259-1262.

\bibitem{BV2}{\sc L.Barbieri-Viale}: Cicli di
codimensione $2$ su variet\`a unirazionali complesse,
preprint (forthcoming with the Ast\'erisque volume
``Proceedings of the $K$-theory
conference'', Strasbourg, 1992), {\it Publ.Dip.Mat.Genova}
{\bf 214}, 1992.

\bibitem{BS}{\sc L. Barbieri-Viale}~and~{\sc V.
Srinivas}: On the N\'eron-Severi group of a
singular variety, {\it J. reine angew. Math.} {\bf 435}
(1993) 65--82.

\bibitem{BSM}{\sc L. Barbieri-Viale}~and~{\sc V.
Srinivas}: The N\'eron-Severi group and the mixed Hodge
structure on $H^2$, to appear on {\it J. reine angew.
Math.}

\bibitem{BBD}{\sc A.Beilinson},~{\sc
J.Bernstein}~and~{\sc P.Deligne}: Faisceaux pervers,
Analyse et topologie sur les espaces singuliers {\it
Asterisque} {\bf 100} (1982).

\bibitem{BE}{\sc A.Beilinson}: Higher regulators and
values of L-functions, {\it J.Soviet Math.} {\bf 30} (1985)
2036-2070.

\bibitem{BMS}{\sc A.Beilinson},~{\sc
R.MacPherson}~and~{\sc R. Schechtman}: Notes on motivic
cohomology, {\it Duke Math.Jour.} {\bf 54} (1987),
679-710.

\bibitem{BL}{\sc S. Bloch}: Lectures on algebraic cycles,
{\it Duke Univ.Math.Ser.} {\bf 4}, Durham, 1980.

\bibitem{BR}{\sc S. Bloch}: Torsion algebraic cycles,
$K_2$ and the Brauer group of function fields, in ``Groupe
de Brauer'', {\it Springer} LNM  {\bf 844} (1981), 75-102.

\bibitem{BO}{\sc S.Bloch~}and{\sc~A.Ogus}: Gersten's
conjecture and the homology of schemes, {\it Ann.Sci.Ecole
Norm.Sup.} {\bf 7} (1974), 181-202.

\bibitem{CO}{\sc
J.-L.Colliot-Th\'el\`ene~}and{\sc~M.~Ojanguren}:
Vari\'et\'es unirationelles non rationelles: au-del\`a de
l'exemple d'Artin et Mumford, {\it Invent.Math.} {\bf 97}
(1989), 141-158.

\bibitem{FU}{\sc Fulton}: Intersection theory, {\it
Springer-Verlag}, Berlin-Heidelberg, 1984.

\bibitem{GIL}{\sc H.Gillet}: Riemann-Roch theorems for
higher algebraic $K$-theory, {\it Adv.in Math.} {\bf 40}
(1981), 203-289.

\bibitem{GID}{\sc H. Gillet}: Deligne homology and
Abel-Jacobi maps, {\it Bull. AMS} {\bf 10}
(1984), 285-288.

\bibitem{GIN}{\sc H.Gillet}: $K$-theory and intersection
theory revisited, {\it K-theory} {\bf 1}
(1987), 405-415.

\bibitem{GO}{\sc R.Godement}: Topologie alg\'ebrique et
th\'eorie des faisceaux, {\it Hermann}, Paris, 1958.

\bibitem{GOM}{\sc M.Goresky}~and~{\sc
R.MacPherson}: La dualit\'e de Poincar\'e pour les espaces
singuliers, {\it C.R.Acad.Sci.Paris} {\bf 284}
S\'erie A (1977), 1549-1551.

\bibitem{GRE}{\sc Grauert~}and{\sc~M.Remmert}: Theory
of Stein spaces, {\it Springer}  Grundlehren Math. Wiss.
{\bf 236}, 1986.

\bibitem{GC}{\sc A.Grothendieck}:~La th\'eorie des
classes de Chern, {\it Bull.Soc.Math.France} {\bf 86}
(1958), 137-154.

\bibitem{GI}{\sc A.Grothendieck}:~Sur quelques
propri\'et\'es fondamentales en th\'eorie des
intersections, in ``S\'eminaire C.Chevalley: Anneaux
des Chow et applications'', E.N.S., 1958.

\bibitem{GR}{\sc A.Grothendieck}:~Le groupe de
Brauer: {\sc i, ii, iii}, in  `` Dix Expos\'es sur la
Cohomologie des Sch\'emas'', {\it North Holland},
Amsterdam, 1968.

\bibitem{HA1}{\sc R.Hartshorne}: Residues and duality
{\it Springer} LNM {\bf 20}, 1990.

\bibitem{HA2}{\sc R.Hartshorne}: On the De Rham
cohomology of algebraic varieties
 {\it Publ. Math.}\, IHES {\bf 29} (1966)
95-103.

\bibitem{HV}{\sc H.Esnault}~and{\sc~E.Viehweg}:
Deligne-Beilinson cohomology, in  Beilinson's
conjectures on special values of $L$-functions, {\it
 Academic Press Perspectives in Math.} {\bf 4} (1987),
43-92.

\bibitem{JAD}{\sc U.Jannsen}: Deligne homology,
Hodge-$\cD$-conjecture and motives, in  Beilinson's
conjectures on special values of $L$-functions, {\it
 Academic Press Perspectives in Math.} {\bf 4}
(1987), 305-372.

\bibitem{JA}{\sc U.Jannsen}: Mixed motives and algebraic
$K$-theory, {\it Springer} LNM {\bf 1400}, 1990.

\bibitem{KS}{\sc M.Kashiwara}~and~{\sc~P.Schapira}: Sheaves
on manifolds, {\it Springer} Grundlehren Math. Wiss. {\bf
292} 1990.

\bibitem{LED}{\sc M.Levine}: Deligne-Beilinson cohomology
for singular varieties, in  Algebraic $K$-theory,
commutative algebra and algebraic geometry, {\it AMS
Contemporary Mathematics} {\bf 126} (1992), 113-146.

\bibitem{Q}{\sc D.Quillen}: Higher algebraic K-theory:~I,
in {\it Springer} LNM  {\bf 341} (1973).

\bibitem{SAM}{\sc P.Samuel}: Relations d'\'equivalence en
g\'eom\'etrie alg\'ebrique, Proceedings ICM, Edinburgh
1958, 470-487.

\bibitem{SH}{\sc C.Sherman}: $\cK$-cohomology of
regular schemes, {\it Comm. Algebra} {\bf 7}
 (1979), 999-1027.

\bibitem{SOU}{\sc C.Soul\'e}: Op\'erations en
$K$-th\'eorie alg\'ebrique, {\it Can.J.Math.} {\bf 37}
(1985), 488-550.

\bibitem{TH}{\sc R.W.Thomason}: Algebraic $K$-theory and
\'etale cohomology, {\it Ann.Sci.Ec.Norm.Sup.} {\bf 18}
(1985), 437-552.

\bibitem{VE}{\sc J.-L.Verdier}: Dualit\'e dans  la
cohomologie des espaces localement compacts {\it
S\'eminaire Bourbaki} {\bf 300}, 1965-66.

 \E{thebibliography}

{}\hfill\\[3cm]

\B{flushright}
\noindent{\sc Universit\`{a} di Genova\\ Dipartimento di
Matematica\\ Via L.B.Alberti, 4\\ Genova 16132 -
Italia\\[4pt] e-mail:\ {\sl barbieri}@{\sl
cartesio.dima.unige.it}}
\E{flushright}

\end{document}